\documentclass[12pt]{article}
\usepackage{amssymb,amsmath}
\usepackage{graphicx}
\usepackage{cite}
\usepackage[british]{babel}

\newcommand{\ud}{\mathrm{d}}
\newtheorem{theorem}{Theorem}[section]
\newtheorem{lemma}{Lemma}[section]

\newtheorem{proposition}{Proposition}[section]
\newcommand{\rem}{\noindent \textbf{Remark. }}
\newcommand{\proof}{\noindent \textbf{Proof. }}
\newcommand{\essinf}{\mathop{{\rm ess~inf}}}
\newcommand{\esssup}{\mathop{{\rm ess~sup}}}
\def\Exp{{\mathbb{E}}}
\def\Pr{{\mathbb{P}}}
\def\R{\mathbb{R}}
\def\Z{\mathbb{Z}}
\def\F{{\mathcal{F}}}
\def\G{{\mathcal{G}}}
\def\A{\mathcal{A}}
\def\S{\mathbb{S}}
\def\SS{{\cal S}}
\def\W{\mathcal{W}}
\def\N{\mathbb{N}}
\def\RR{\mathcal{R}}
\def\eps{{\varepsilon}}
\def\1{{\bf 1}}
\def\bx{{\bf x}}
\def\by{{\bf y}}
\def\bu{{\bf u}}
\def\bv{{\bf v}}
\def\be{{\bf e}}
\def\0{{\bf 0}}
\def\sgn{{\rm sgn }}

 \newcommand{\qed}{\hfill $\Box$}

\numberwithin{equation}{section}

\allowdisplaybreaks

\author{Iain M. MacPhee$^1$ \and Mikhail V. Menshikov\footnote{Corresponding author; email \texttt{Mikhail.Menshikov@durham.ac.uk}}$^{~,1}$
\and
Andrew R. Wade$^{2}$}

\title{Angular asymptotics for multi-dimensional non-homogeneous random walks with asymptotically zero drift}

\begin{document}

\maketitle

{\footnotesize
\noindent
$^1$
Department of Mathematical Sciences,
University of Durham,
South Road, Durham DH1 3LE, UK.\\
\noindent
$^2$
Department of Mathematics and Statistics,
 University of Strathclyde,
 26 Richmond Street, Glasgow G1 1XH, UK. }

\begin{abstract}
We study the first exit time $\tau$ from an
arbitrary cone with apex at the origin
by a non-homogeneous
random walk (Markov chain)
on $\Z^d$ ($d \geq 2$)
with mean drift that is asymptotically zero.
Specifically, if the mean drift at $\bx \in \Z^d$ is of magnitude
$O(\| \bx\|^{-1})$,
we show that
$\tau<\infty$ a.s.~for any cone. On the other hand,
for an appropriate drift field with mean
drifts of magnitude $\| \bx\|^{-\beta}$, $\beta \in (0,1)$,
we prove that our random walk has a limiting (random) direction and so
eventually remains in an arbitrarily narrow cone.
The conditions imposed on the random walk are minimal: we assume
only a uniform bound on $2$nd moments for the increments and
a form of weak isotropy. We give several illustrative examples,
including  a random walk in random environment model.
\end{abstract}

\vskip 2mm


\vskip 2mm


\section{Introduction}

The theory of time- and space-homogeneous random walks on $\Z^d$ ($d
\geq 2$), i.e., sums of i.i.d.~random integer-component vectors, is
classical and
extensive; see for example \cite{spitzerrw,cohenbook,lawlerbook}.  For
random walks that are not spatially homogeneous the theory is less
complete, and many techniques available for the study of homogeneous
random walks can no longer be applied, or are considerably
complicated; see, for instance, \cite{lawler91,mustapha}. In the present
paper we study angular properties of non-homogeneous
random walks, specifically exit times from cones and existence
of limiting directions.

In general non-homogeneous processes can be wild;
thus we restrict ourselves to walks that have
mean
drift that tends to zero as the distance to the origin tends
to infinity (but with no restriction on the direction of the drift)
and satisfy some  weak regularity conditions on the jumps.
We do not impose on the increments of the random walk conditions of boundedness, symmetry, or uniform ellipticity,
as are assumed, for example, for the results on non-homogeneous
random walks in \cite{lawler91,mustapha}. Importantly, we do not impose any direct
restrictions on the correlation structure of the components
of the increments of the process.
Random walk models are applied in many contexts. Often, simplifying
assumptions of homogeneity are made in order to make such models
tractable, whereas non-homogeneity
is more realistic. Thus our non-homogeneous
model shares some motivation
 with random walks in random environments (see e.g.~\cite{zeitouni}); in such terms,
 our results
 deal with a particular class of `asymptotically zero drift'
 environments (cf.~Example 5 in  Section \ref{examples} below).
  In the present
paper we develop methods to study passage-times for certain sets
for such non-homogeneous random walks.

We now describe informally the
type of non-homogeneous random walk studied in the present paper.
Consider a Markovian
random walk on $\Z^d$, homogeneous in time but not necessarily in space,
so that the transition function depends upon the walk's current
location.
Suppose that the walk
 has one-step
 mean drift function that
tends to zero as the
distance from the origin
tends to infinity. This asymptotically
zero drift regime  is the natural
setting in which
  to probe the transition
away from behaviour that is
essentially `zero-drift'
in character. In one dimension,
the corresponding
regime
is
 rather well understood,
 following fundamental work of Lamperti; see
\cite{lamp1,lamp2,mai,mp,mvw} and the Appendix in \cite{aim}
(some analogous results in the continuous setting of
 Brownian motion with asymptotically zero drift
are given more recently in \cite{dante2}).
 Problems
in higher dimensions of
a `radial' nature
can often be reduced to this one-dimensional
case. The exit-from-cones problems that we consider
in the present paper (which we describe below), on the other hand,
are to a large extent `transverse' (and inhomogeneous)
in nature and so
are truly
many-dimensional.
Moreover, the many-dimensional
case is qualitatively different from the one-dimensional case
(see Theorem \ref{thm1} below).

The random walks that we consider are non-homogeneous, but
some regularity assumptions are certainly  required
for  our results.
We assume a {\em weak isotropy} condition
without which highly degenerate behaviour is possible.
In addition,
we restrict our attention to random walks on unbounded subsets of
 $\Z^d$ with
some moment condition on the jumps.
 We need some
regularity conditions on the state-space of our walk
and   it is most convenient to take the structure of $\Z^d$. We are confident that
our proofs can be adapted for more general state spaces.

Our main theorems can be summarized as follows: (i) a walk
with mean drift of magnitude $O(\| \bx \|^{-1})$
at $\bx$ will leave any cone in finite time almost surely
(and indeed hit any cone), while (ii)
an appropriate drift field with
magnitude of order $O(\| \bx \|^{-\beta})$, $\beta \in (0,1)$
can lead to the existence of an asymptotic direction for the walk
(so that it eventually remains in an arbitrarily thin cone). Note that
the class of random walks with mean drifts $O(\| \bx \|^{-1})$
to which result (i) applies is very wide: such a walk can be
 transient, null-recurrent, or positive-recurrent
 (cf.~\cite{lamp1,lamp2}) and can be diffusive or sub-diffusive
 (cf.~\cite[Section 4]{mvw}.)

Before stating our theorems formally,
we briefly describe some of the relevant existing literature.
The theory of homogeneous zero-mean random walks stands
hand-in-hand
with  the corresponding
continuum theory for Brownian motion. Once
the assumption of spatial homogeneity
is removed,
Brownian motion ceases to
be a reliable analogy for the random
walk problem. In the case of one
dimension, this is exemplified
by results
on processes with asymptotically
zero mean drifts; see e.g.~\cite{lamp1,lamp2}.
 For the non-homogeneous random walks considered in the present paper, we will demonstrate
behaviour substantially different to that of standard Brownian motion.

In \cite{wedges} the authors give conditions under which our non-homogeneous random
walk does display essentially `Brownian' behaviour.
The study of the exit-time of standard Brownian motion from cones
goes back at least to Spitzer \cite{spitzer} and a deep
analysis was undertaken by Burkholder \cite{burkh}; see \cite{bansmi}
for some more recent work.
The random walk case has received less attention.
A  body of work by Varopoulos starting with \cite{var1} deals
with exit-from-cones problems for random walks that
have {\em mean drift zero} but are (at least for
some of the results in \cite{var1}) allowed
to be non-homogeneous. In \cite{var1}, finer behaviour (such as tails
of exit times) was studied, and consequently the conditions on the walks
imposed in \cite{var1} are stronger than ours in several respects, such as
an assumption of orthogonality on the covariance structure of the increments.

In the next section we
give the precise formulation of the model, our main results, and a discussion.
In particular, in Section \ref{descr} we
formally define our model and our assumptions.
In Section \ref{results} we state our main results.
Then in Section \ref{examples} we give several examples of processes to which our theorems can be applied,
including `centrally biased random walks', half-plane excursions,
and a random walk in random environment model.
In Section \ref{open} we mention some possible directions for future research.
Finally, in Section \ref{outline} we give a brief outline of the
technical part of the paper, which contains the proofs of our
results.

\section{Model, results, and discussion}
\label{model}

\subsection{Description of the model}
\label{descr}

In this section we describe
more precisely
the probabilistic model that is our object of study.
First
we collect some notation.
Throughout we assume $d \in \{2,3,\ldots \}$ and
work in $\R^d$; $2$ is the minimum number of dimensions
in which the phenomena that we study appear, although
analogues of our results in the case $d=1$
are in a sense provided by Lamperti \cite{lamp1,lamp2}.
For $\bx \in \R^d$,  write
$\bx = (x_1, \ldots,x_d)$ in Cartesian coordinates.
 Let $\| \cdot \|$   denote the Euclidean
norm on $\R^d$. For a non-zero
vector $\bx \in \R^d$ we use the usual notation $\hat \bx := \bx / \| \bx \|$ for the
corresponding unit vector. Write $\0 := (0,\ldots,0)$
for the origin
and $\be_1,\ldots, \be_d$ for the standard orthonormal
basis of $\R^d$. For vectors $\bu, \bv \in \R^d$
we use $\bu \cdot \bv$ to denote
their scalar product.

Let $\Xi = (\xi_t)_{t \in \Z^+}$ be a discrete-time
Markov
 process
 with state-space $\SS$ an unbounded subset of
 $\Z^d$. Since we are concerned
 crucially with the spatial aspects of the process,
 it is natural to view our process a {\em random walk} on $\SS \subseteq \Z^d$,
 although it will certainly not, in general, be a sum
 of i.i.d.~random vectors.
The random walk $\Xi$ will be time-homogeneous
but not necessarily space-homogeneous;
we will impose some natural regularity
assumptions on
the increment distribution for our walk, which we describe next.

We need to impose some form of regularity condition that
ensures the walk cannot become trapped in lower-dimensional
subspaces  or finite sets. To this end,
we will assume the following  weak isotropy condition:

\begin{itemize}
\item[(A1)] There exist $\kappa>0$, $k  \in \N$ and
  $n_0 \in \N$ such that
\[ \min_{ \bx \in \SS} \min_{\by \in \{ \pm k  \be_i , ~i=1,\ldots, d \} }
\Pr [ \xi_{t+n_0} - \xi_t = \by \mid \xi_t = \bx ] \geq \kappa ~~~(t \in \Z^+).\]
\end{itemize}

Note that (A1) is an $n_0$-step regularity
condition. In terms of one-step
regularity, its implications are
minimal:
a simple consequence of (A1) is that for any $\bx \in \SS$
\[ \Pr [ \xi_{t+1} = \bx \mid \xi_t = \bx ]
= ( \Pr[ \xi_{t+n_0} = \bx, \ldots, \xi_{t+1} = \bx \mid \xi_t = \bx]
)^{1/n_0} \leq (1-2d\kappa)^{1/n_0} \leq 1 - (2d\kappa/n_0) \]
(note $\kappa \leq 1/(2d)$) so that
\[
\Pr [ \xi_{t+1} \neq \bx \mid \xi_t =\bx ] \geq (2d\kappa/n_0) > 0 \]
uniformly in $\bx$. Condition (A1) can be seen as a form of ellipticity,
but is weaker than uniform ellipticity (such as often
assumed in the random walk in random environment literature, see
e.g.~\cite{zeitouni}). For example, there can be sites $\bx \in \SS$
at which the jump distribution degenerates completely and the walk moves
deterministically. (When later we discuss walks with asymptotically zero mean
drift, for all $\bx$ with $\| \bx\|$ large enough this extreme degeneracy
is excluded, although the jump distribution at $\bx$ may still
be supported on a lower-dimensional subspace.) At first sight it seems
that we are
losing some generality in (A1) by enforcing a single $k$ and $n_0$
for each of the $2d$ directions in the condition --- but this
is not in fact any sacrifice, as we show in Proposition \ref{propA1}
below (see Section \ref{open}). Finally, note that another
consequence of (A1) is that $\limsup_{t \to \infty} \| \xi_t \| = \infty$ a.s..

Our time-homogeneity and Markov assumptions imply that
the distribution of the increment
$\xi_{t+1} - \xi_t$ depends only on the position $\xi_t$ and not $t$.
Our second regularity condition is an assumption of
finiteness of second moments for the increments
of $\Xi$:
\begin{itemize}
\item[(A2)] There exists $B_0 \in (0,\infty)$ such that
\[  \max_{\bx \in \SS} \Exp [ \| \xi_{t+1} - \xi_t \|^{2} \mid \xi_t = \bx ] \leq B_0 .\]
 \end{itemize}
 It is interesting that for our theorems and
 with our techniques $2$ moments suffice, rather than $2+\eps$
 moments or uniformly bounded jumps as are often assumed in similar situations.
 Under (A2), the mean of $\xi_{t+1} - \xi_t$ given $\{ \xi_t = \bx\}$ is well-defined.
Denote the one-step
{\em mean drift} vector
$\mu (\bx) := \Exp [ \xi_{t+1} - \xi_t \mid \xi_t = \bx ]$
for $\bx \in \SS$.
  We are primarily
interested in the case where
 the
random walk
 has {\em asymptotically zero} mean drift, i.e.,
$\lim_{\| \bx \| \to \infty} \| \mu (\bx)\| =0$.

Write $\S_d := \{ \bu \in \R^d : \| \bu \| = 1 \}$ for the unit
sphere in $\R^d$.
For $\bu \in \S_d$,
$\alpha \in (0,\pi)$
let $\W_d(\bu; \alpha)$ be an open (circular)
cone in $\R^d$ with apex $\0$, principal direction $\bu$,
and half-angle $\alpha$:
\[ \W_d (\bu ; \alpha) := \{ \bx \in \R^d : \bu \cdot \hat \bx >  \cos \alpha \}.\]

A central quantity in this paper is
the random walk's first exit time from the cone $\W_d(\bu;\alpha)$ (starting
from inside the cone).
Define the random time
\begin{align*}
 \tau_\alpha
 := \min \{ t \in \Z^+ : \xi_t \notin  \W_d (\bu ; \alpha )
    \}.\end{align*}
The notation $\tau_\alpha$ suppresses the dependence
on the starting point $\xi_0$ and the cone direction $\bu$.
 Note that the complementary cone $\R^d \setminus \W_d ( \bu ; \alpha)$
has interior $ \W_d ( -\bu ; \pi-\alpha)$ so exit from a large cone is equivalent to hitting
a small cone. Exit from a small cone does not in general imply hitting {\em any} small
cone for a non-homogeneous random walk without some condition that prevents confinement of the
walk to a subspace of $\R^d$. This is why we need a condition such as (A1).\\

\rem The time-homogeneity and Markov assumptions that we make are not crucial
for our results, and are not essentially used in our proofs. However, to avoid
complicating the statements of our theorems we have not used the
maximum generality in this respect. In fact, we essentially prove our Theorem
\ref{thm2} in the more general setting (see Section \ref{sectran}).

\subsection{Main results}
\label{results}

Our first result, Theorem \ref{thm1} below,
deals with the  case where the mean drift is $O(\| \bx\|^{-1})$; we will see
that this case is critical for our properties of interest.
It is often useful to view our general model as a perturbation of the
zero-drift case.
It is perhaps intuitively clear,
by analogy with Brownian motion, that a zero-drift homogeneous
random walk
 on $\Z^d$ satisfying  suitable regularity conditions will
 exit any cone
in almost surely finite time. Note that care is needed even in the zero-drift case,
since random walks with zero drift can behave very differently from Brownian motion
due to correlation structure of the increments: see e.g.~\cite{klein}.
 It is less clear that such a result is true
for random
walks that are non-homogeneous and have an arbitrary correlation
structure for their increments.   Theorem \ref{thm1}  provides the much stronger result that
 the exit time is a.s.~finite in the {\em asymptotically} zero drift
setting provided that the mean drift is
$O(\| \bx\|^{-1})$. Moreover, if this latter condition fails,
the result may be false (see Theorem \ref{thm2} below); in this
sense, Theorem \ref{thm1} is best possible.

\begin{theorem}
\label{thm1}
Suppose that (A1) and (A2)
 hold, and
that for $\bx \in \SS$ as $\|\bx\| \to \infty$,
\begin{align}
\label{drift1}
\| \mu (\bx)\| = O ( \| \bx \|^{-1} ) .\end{align}
Then for any $\alpha \in (0,\pi)$, any $\bu \in \S_d$,
 and any $\bx \in \SS \cap \W_d (\bu ; \alpha)$
\[ \Pr [ \tau_\alpha < \infty \mid  \xi_0 = \bx ]  = 1. \]
 \end{theorem}

As a special case, Theorem \ref{thm1} includes the case of a non-homogeneous
random walk with zero drift.
The only similar result that we could find explicitly stated
in the literature is
in \cite{var1}, where it was shown that $\tau_\alpha < \infty$ a.s.~for
a non-homogeneous random walk with mean drift zero under a
condition of uniformly bounded
jumps and several other technical conditions including assumptions on correlation
structure of the jumps and conditions on the reversed process.
Thus Theorem \ref{thm1} provides a proof of the result $\tau_\alpha < \infty$ a.s.~in the
zero drift setting under
conditions that are weaker in several directions (in particular
the assumptions on the increments) than those in \cite{var1}.
The main object
of \cite{var1} was to address the more delicate question of obtaining tight bounds
for the tail of $\tau_\alpha$. In our more general setting  (with mean drift
{\em asymptotically} zero) the tails of $\tau_\alpha$ depend
crucially on the drift field, even in the case where the mean drift is
$O(\| \bx\|^{-1})$, or indeed identically zero, and we do not consider the problem of tail bounds in the present
paper. However, in \cite{wedges} the authors do show that, for $d=2$,
if the mean drift is $O(\| \bx\|^{-1})$
then $\tau_\alpha$ has a {\em polynomial} tail under
the condition of uniformly bounded jumps. For an informative example of the impact of
increment correlation structure on the existence of
exit-time moments in a simple setting, see \cite{klein}.

We emphasize that walks $\Xi$ satisfying Theorem \ref{thm1}
can display a wide range of behaviour. For example, in the case
of radial drift $\mu (\bx) = c \| \bx \|^{-1} \hat \bx$ for $c \in \R$,
it can be shown by an analysis of the process $\| \xi_t \|$
(possibly under some additional regularity assumptions) that, depending on $c$,
$\Xi$ can be positive-recurrent, null-recurrent, or transient
 (see Example 2 in Section \ref{examples} below)
and that $\Xi$ can be diffusive or sub-diffusive (see
\cite[Section 4]{mvw}).

 Theorem \ref{thm1} contrasts
sharply with the situation in one dimension
\cite{lamp1,lamp2}, where a drift of $O(x^{-1})$ at $x$
does not imply  finiteness of the time
of exit from a half-line. In $d=2$,
 Theorem \ref{thm1}
 gives information on the winding of the walk around the origin; an early result
on the
winding number of planar Brownian motion
is also contained in Spitzer's paper \cite{spitzer} and
a more recent reference,
including
corresponding
results for homogeneous random walks,
is \cite{shi}. Theorem \ref{thm1} generalizes such winding properties
naturally to higher dimensions.

Now we move on to the supercritical case.
Theorem \ref{thm2} shows that for a radial drift
field, with outwards drift
greater in order than $\| \bx \|^{-1}$,
the walk now has a limiting direction,
in complete contrast to the situation in Theorem
\ref{thm1}. In other words, the random walk
eventually remains in an arbitrarily thin cone.

 \begin{theorem}
 \label{thm2}
 Suppose that (A1) and (A2) hold.
 Suppose
that for some $\beta \in (0,1)$, $c>0$, $\delta>0$,
 and $A_0>0$,
 \begin{align}
 \label{trancona}
  \min_{\bx \in \SS : \| \bx \| > A_0 } \{ \| \bx \|^\beta \mu(\bx) \cdot \hat \bx \} \geq c ,
 ~\textrm{and}\\
 \label{trancon2a}
 \max_{\bx \in \SS, \| \bx\| > A_0 } \sup_{ \bu \in \S_d : \bu \cdot \bx = 0}
   \| \bx \|^{\beta+\delta} | \mu(\bx) \cdot \bu | < \infty .
 \end{align}
 Then for any $\xi_0 \in \Z^d$ we have that $\| \xi_t \| \to \infty$ a.s., and
 there exists a random unit vector $\bv \in \S_d$,
 whose distribution is supported on all of $\S_d$,
  such that  a.s.~as $t \to \infty$
 \[ \frac{\xi_t}{\| \xi_t \|} \to \bv. \]
 \end{theorem}

Note that Theorem \ref{thm2} says that the random walk is transient,
a fact that does not follow immediately from known results
(for instance to apply Lamperti's results \cite{lamp1}
one needs a stronger moment assumption than (A2)).
 A natural example to which Theorem \ref{thm2}
 applies is a walk where for some $c>0$ and $\beta \in (0,1)$
 \[ \mu (\bx) \cdot \hat \bx = c \| \bx \|^{-\beta} , ~~~
 {\rm and}~~~ | \mu (\bx) \cdot \bu | = O ( \| \bx \|^{-1} ) ,\]
 for all $\bu$ orthogonal to $\bx$. See also Example 2 in Section \ref{examples} below.

 Theorem \ref{thm2} covers walks that are sub-ballistic
  (i.e.~have zero speed, asymptotically). We could not
 find results on limiting directions for non-homogeneous
 random walks in the literature. The phenomenon
 of limiting direction for homogeneous walks on spaces more exotic
 than $\Z^d$ has been studied: see e.g.~\cite{karlsson} and
 references therein.
  In the next section we illustrate our two main results with
  some examples.

 \subsection{Examples and comments}
 \label{examples}

We now list some particular examples
of random walks with which we
will illustrate Theorems \ref{thm1} and \ref{thm2}.
In some cases we  assume the following slightly
stronger version of (A2):
\begin{itemize}
\item[(A2+)] There exist $\eps >0$ and $B_0 \in (0,\infty)$ such that
\[  \max_{\bx \in \SS} \Exp [ \| \xi_{t+1} - \xi_t \|^{2+\eps} \mid \xi_t = \bx ] \leq B_0 .\]
 \end{itemize}

\subsubsection*{Example 1: Zero-drift non-homogeneous random walk.}

Let $d \geq 2$.
Suppose
that (A1) and (A2) hold and
$\mu(\bx) \equiv \0$. Note that even for this
example, the random walk is not necessarily homogeneous
and the covariance structure of the increments
is arbitrary,
so the walk is not covered by classical work such
as \cite{spitzerrw} or more recent work such as
\cite{var1}. One can construct examples of such walks that are transient in $d=2$, or recurrent in $d \geq 3$, for instance.
Theorem \ref{thm1} immediately implies that
in this case the walk leaves
any cone in finite time.

\subsubsection*{Example 2:
Random walk with radial drift.}

\begin{figure}
\centering
\begin{minipage}[c]{0.58\textwidth}
\centering \includegraphics[width=7cm]{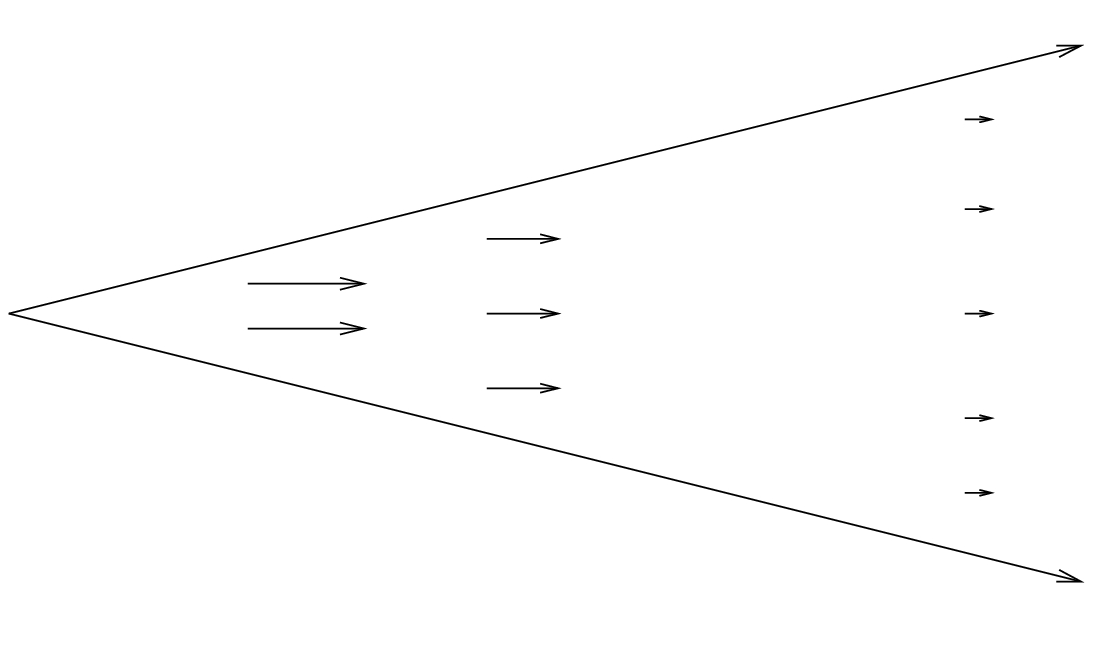}
\end{minipage}
\quad
\begin{minipage}[c]{0.38\textwidth}
\centering \includegraphics[width=5cm]{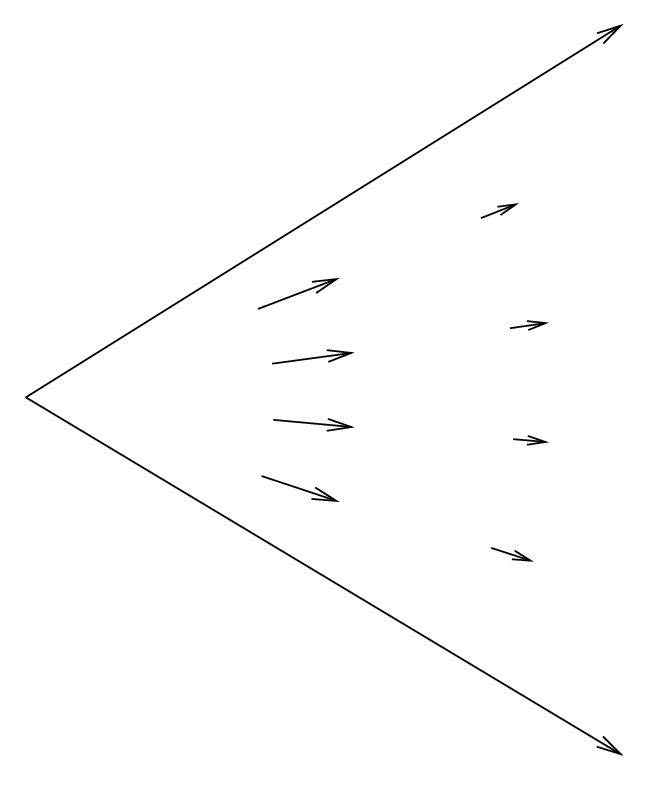}
\end{minipage}
\caption{Two examples of drift fields: $c x_1^{-1} \be_1$ (left) and $c \| \bx \|^{-1} \hat \bx$ (right).}
\label{fig1}
\end{figure}

Let $d \geq 2$. Suppose that (A1) and (A2)  hold,
and that for some $c \in \R$, $\beta>0$,
for $\bx \neq \0$, $\mu(\bx) = c \| \bx\|^{-\beta} \hat \bx$.
An example of a suitable drift field (for $d=2$) is illustrated
in the second part of Figure \ref{fig1}.
This kind of model has
been called a {\em centrally biased} random walk
(see e.g.~\cite[Section 4]{lamp1}). The following result is again immediate
from Theorems \ref{thm1} and \ref{thm2}.

\begin{theorem}
\label{thm3}
Suppose $\Xi$ is as in Example 2. Let $\alpha \in (0,\pi)$ and $\bu \in \S_d$.
\begin{itemize}
\item[(i)] If $\beta \geq 1$, then for  any $\bx \in \W_d(\bu; \alpha)$,
$\Pr [ \tau_\alpha < \infty \mid \xi_0 = \bx ] = 1$.
    \item[(ii)] If $\beta < 1$ and $c>0$, then for any $\bx \in \W_d(\bu ; \alpha)$,
    $\|\xi_t\| \to \infty$ a.s.~and
    $\xi_t/\| \xi_t \| \to \bv$ a.s.~as $t \to \infty$,
    for some $\bv$ with distribution supported on $\S_d$.
        \end{itemize}
    \end{theorem}

It is worth comparing the behaviour of the walk in this example
in terms
of exit from cones to its  recurrence/transience
behaviour (in terms of returning to bounded sets), which can be
obtained from study of the process $\| \xi_t \|$.
Results of Lamperti \cite{lamp1,lamp2} (see also \cite{aim,mai})
imply that, at least if we assume (A2+),
\begin{itemize}
\item If $\beta > 1$, $\Xi$ is recurrent in $d=2$ and transient for $d \geq 3$;
\item If $\beta < 1$, then $\Xi$ is transient for $c>0$ and positive-recurrent for $c<0$.
\end{itemize}
The case $\beta=1$ is critical from the point of view of the
recurrence classification
 (see in particular
the discussion
around (4.13) in \cite{lamp1}),
and, for any $d$, $\Xi$ can be either positive-recurrent, null-recurrent, or transient,
depending on $c$. In particular, there exist $c_0, c_1 \in (0,\infty)$ (depending
on $d$ and $B_0$) such that
  $\Xi$   is positive-recurrent for $c<-c_1<0$   but
   transient
  for $c>c_0>0$.
  Thus when $\beta=1$ and $c>c_0$,
  $\Xi$ is transient
 and so eventually
leaves every bounded region, but,
on the other hand (by Theorem \ref{thm3}(i))
such a walk
will also eventually  leave any wedge.
In other words, although $\|\xi_t\| \to \infty$
the walk has no limiting direction.

\subsubsection*{Example 3:
 Random walks with drift
in the principal direction.}

Let $d \geq 2$.
It is interesting to contrast two
apparently similar types of random walk
on the half-space $\W_d(\be_1;\pi/2)$. Suppose that (A1) and
(A2+) hold. Suppose either
\begin{itemize}
\item[(a)]
for some $c \in \R, \beta>0$, for $\bx \in \W_d (\be_1 ; \pi/2)$,
$\mu (\bx) = c \| \bx\|^{-\beta} \be_1$; or
\item[(b)]
for some $c \in \R, \beta>0$, for $\bx \in \W_d (\be_1 ; \pi/2)$,
$\mu (\bx) = c x_1^{-\beta} \be_1$.
\end{itemize}
An example of a suitable drift field in case (b) ($d=2$) is illustrated
in the first part of Figure \ref{fig1}.  The following result
  is again a consequence of Theorems \ref{thm1} and \ref{thm2}, but requires
  some extra work: we present its proof at the end of Section \ref{prelim}.

  \begin{theorem}
\label{thm4}
Suppose $\Xi$ is as in Example 3. Suppose $\beta =1$.
If $\alpha \in (0,\pi/2)$, then in either case (a) or (b), for
 any $\bx \in \W_d(\be_1 ; \alpha)$,
\[ \Pr [ \tau_\alpha < \infty \mid \xi_0 = \bx ] = 1.\]
On the other hand,
in case (a), for
 any $\bx \in \W_d(\be_1;\pi/2)$,
\[ \Pr [ \tau_{\pi/2} < \infty \mid \xi_0 = \bx ] = 1,\]
but in case (b) there exists $c_0 \in (0,\infty)$
such that for $c>c_0$ and any $\bx \in \W_d(\be_1;\pi/2)$,
\[ \| \xi_t \| \to \infty, ~~~{\rm and }~~~  \xi_t \cdot \be_1 \to \infty ~~~{\rm a.s.}.
\]
    \end{theorem}

Theorem \ref{thm4} shows that the
difference in qualitative
behaviour between cases (a) and (b)
is manifest in terms of
leaving the  half-space.
In particular, when the mean drift
is $c/x_1$ in the $\be_1$ direction, the walk leaves a wedge of angle
$\alpha < \pi/2$, but, for $c$ large enough, with positive probability eventually remains in
the half-plane.
However when the mean drift
is $c/\| \bx \|$ in the $\be_1$ direction, the
 walk always leaves the wedge,
  even when $\alpha = \pi/2$.
The instance of Example 3, case (b),
when $\alpha = \pi/2$
demonstrates homogeneity
in the $\be_2$ direction, and so
is  related to
the one-dimensional so-called Lamperti problem
named after \cite{lamp1,lamp2}. In the case $\alpha=\pi/2$, case (a) demonstrates
a more
localized perturbation, since near
the boundary of the half-plane we can have $\| \bx \| \gg x_1$.

The primary interest of Example 3 is the case $\beta=1$.
For reasons of space we do not consider here the case $\beta \in (0,1)$ of Example 3 (either
(a) or (b)); we expect that this case too can be studied using our methods.

\subsubsection*{Example 4: Random walk half-plane excursion.}

We point out a particularly simple case of Example 3, case (b) above,
which is of interest in its own right. This is
the so-called random walk half-plane excursion
(see \cite{lawlerconf}, pp.~1--2). This process
is obtained, loosely speaking, by conditioning
 a simple symmetric random walk on $\Z^2$
never to exit a half-plane: see \cite{lawlerconf} for details.
The construction readily extends to general dimensions $d \geq 2$,
but for simplicity we discuss the planar case.
In this case $\Xi$ has transition probabilities
\begin{align*}
 \Pr [ \xi_{t+1} = ( x_1, x_2 \pm 1 ) \mid \xi_t = (x_1,x_2) ] & = \frac{1}{4} \\
 \Pr [ \xi_{t+1} = (x_1 \pm 1 , x_2 ) \mid \xi_t = (x_1,x_2) ] & = \frac{x_1 \pm 1}{4x_1} ,
 \end{align*}
 for $(x_1,x_2) \in \Z^2$, $x_1 \geq 1$. Hence
 \[ \mu (\bx) = \frac{1}{2 x_1} \be_1 ,\]
 and we are in the case of Example 3(b) as described above.
 Theorem
 \ref{thm4} implies that for any $\alpha \in (0,\pi/2)$,
 the walk leaves the wedge $\W_2(\be_1; \alpha)$ in finite time almost surely.
 On the other hand, note that $\Xi$ is transient and in fact $\xi_t \cdot \be_1 \to \infty$
 almost surely,  by for instance
 Lamperti's results \cite{lamp1} (in fact one can take $c_0=1/4$ in Theorem \ref{thm4} above,
 so the final statement of that theorem applies: see the proof of Theorem \ref{thm4}
 in Section \ref{prelim}).

\subsubsection*{Example 5: Random walk in   random environment.}

We give a final example of a slightly different flavour.
Let $d \geq 2$. Suppose that
each site $\bx \in \Z^d$ carries random
$d$-vectors ${\bf Y}^\bx$ and $\chi^\bx$, all independent,   where
${\bf Y}^\bx = (Y^\bx_1, \ldots, Y^\bx_d)$ has an arbitrary
distribution (possibly even dependent on $\bx$)
on the simplex $\{ (y_1, \ldots, y_d) \in [0,\infty)^d : y_1 + \cdots + y_d = 1 \}$,
and $\chi^\bx = (\chi^\bx_1,\ldots,\chi^\bx_d)$ is an independent copy of $\chi = (\chi_1,\ldots,\chi_d)$, whose components $| \chi_i |$
are bounded uniformly in $i$.
Let $\omega := ( ({\bf Y}^\bx , \chi^\bx ) )_{\bx \in \Z^d}$
be the {\em random environment}. Given $\omega$, define a
nearest-neighbour Markovian random walk $\Xi$ on $\Z^d$ via its
transition law $\Pr_\omega$ given by, for $i \in \{1,\ldots,d\}$,
\[ \Pr_\omega [ \xi_{t+1} - \xi_t = \be_i \mid \xi_t = \bx ] =
\frac{1}{4d} + \frac{Y^\bx_i}{4} + \frac{\chi^\bx_i}{\| \bx \| } ,
~~~ \Pr_\omega [ \xi_{t+1} - \xi_t = -\be_i \mid \xi_t = \bx ] =
\frac{1}{4d} + \frac{Y^\bx_i}{4} - \frac{\chi^\bx_i}{\| \bx \| } , \]
unless either of these quantities lies outside the interval $[\frac{1}{4d}, 1 - \frac{1}{4d}]$,
in which case we replace both probabilities in question by $\frac{1}{4d} + \frac{Y^\bx_i}{4}$ (for
almost every $\omega$,
this modification
will only  apply within a finite ball around the origin).
 Thus $\Xi$ is a random walk in random environment (RWRE).
Then, given $\omega$, $\mu(\bx) = 2 \| \bx \|^{-1} ( \chi^\bx_1, \ldots, \chi^\bx_d )$,
so that $\| \mu (\bx) \| = O (\| \bx \|^{-1})$, uniformly for almost every $\omega$,
by the conditions on $\chi$. Thus a consequence of Theorem \ref{thm1}
is that for almost every $\omega$, for any $\alpha \in (0,\pi)$,
$\tau_\alpha < \infty$ a.s.. To the best of the authors' knowledge, the recurrence/transience
classification of this RWRE is at present
 an open problem. An analogous model
in $d=1$ where $Y^\bx_1 \equiv 1$ for all $\bx$
(random perturbation of the simple symmetric random walk
on $\Z^+$) was studied in \cite{mw1}: Theorem 2
parts (iii)--(v) in \cite{mw1}
give the complete recurrence classification in that case.

\subsection{Extensions, open problems, and further remarks}
\label{open}

As we have already indicated, we essentially prove Theorem \ref{thm2}
without the assumptions of time-homogeneity or the Markov
property (see Section \ref{sectran} below). It should be
possible to prove an appropriate extension
of Theorem \ref{thm1} in similar generality. The assumption of the state-space
being $\Z^d$ is not essentially used in the proof of Theorem \ref{thm2},
which we could have stated for more general walks on $\R^d$ under
an appropriate analogue of (A1); the state-space assumption is
central to the decomposition idea in the proof of Theorem
\ref{thm1} (see Section \ref{prf1}), but we believe that the method
should extend to more general state-spaces
assuming an appropriate generalization of the isotropy
condition (A1).

 As  mentioned above, condition (A1) is more general than it might
 first appear. In fact it is equivalent to the following.

\begin{itemize}
\item[(A1$^\prime$)] There exist $\kappa>0$, and $k_i \in \N$, $n_i \in \N$ for
$i \in \{ \pm 1, \pm 2,\ldots, \pm d\}$ such that
\[ \min_{ \bx \in \SS} \min_{i \in \{ \pm 1, \ldots, \pm d \} }
\Pr [ \xi_{t+n_i} - \xi_t = k_i \sgn (i) \be_{|i|} \mid \xi_t = \bx ] \geq \kappa ~~~(t \in \Z^+).\]
\end{itemize}

 \begin{proposition}
 \label{propA1}
 Conditions (A1) and (A1$^\prime$) are equivalent.
 \end{proposition}

 We prove Proposition \ref{propA1} in Section \ref{prelim}.
It seems
unlikely that the conditions (A1) and (A2)
can be relaxed
to any significant degree
 in Theorems \ref{thm1} and \ref{thm2}.
If (A1) is absent,  Theorem
\ref{thm1} may fail by the walk getting trapped in
a low-dimensional subspace.
For example, if the only possible jumps
of the walk are in the $\pm \be_1$ directions, it will be trapped
on a line. Then one-dimensional results (see e.g.~\cite{lamp1})
imply that even for a mean drift of magnitude $O(\|\bx\|^{-1})$
the process can be transient in the positive $\be_1$ direction, and
so will with positive probability never leave any cone with principal
axis in the $\be_1$ direction, contradiction Theorem \ref{thm1}.

In Theorem \ref{thm2}, some condition such as (A1)
is needed to ensure that $\limsup_{t \to \infty} \| \xi_t \| = +\infty$ a.s.,
or else the walk can get stuck in a finite ball around the origin before
the drift asymptotics take effect. We suspect that the moment
condition (A2) is close to optimal in Theorem \ref{thm1}.
It seems likely that in Theorem \ref{thm2}, (A2) can be replaced
by a uniform bound on $1+ \beta + \eps <2$ moments $(\eps>0)$,
by a more delicate analysis in Lemma \ref{hsup} below. To avoid additional
complications, here we are satisfied with the uniform assumption (A2) throughout.

Several open problems remain. Perhaps the most
interesting, and the natural next question to address,
is the study of the tails
(or moments)
of $\tau_\alpha$
when $\| \mu (\bx)\| = O (\| \bx \|^{-1})$. It is not hard to see
(for instance by comparison with one-dimensional results such as \cite{lamp2,aim})
 that there exists
a wide array
of possible tail behaviours for $\tau_\alpha$.
The authors have studied the case $d=2$, under some additional
assumptions, in \cite{wedges}: of course, covariance
structure of the increments is crucial here
(cf.~\cite{klein}). In particular,
in \cite{wedges} we show that in $d=2$, when
$\| \mu (\bx) \| = o (\| \bx \|^{-1})$
the tails of $\tau_\alpha$ are, to first order,
the same as in the Brownian motion case
under assumptions on correlations
(cf.~Spitzer's theorem \cite{spitzer}). However, the
general picture when $\| \mu (\bx) \| = O ( \| \bx \|^{-1})$ is far from
complete
even in $d=2$.

\subsection{Paper outline}
\label{outline}

The outline of the remainder of the paper is as follows.
In Section \ref{prelim} we collect some preparatory results
and prove Proposition \ref{propA1} and
Theorem  \ref{thm4}. Sections \ref{prf1}
and \ref{sectran}
are devoted to the proofs of Theorems \ref{thm1} and \ref{thm2}
respectively. The two proofs are essentially independent,
so either of these two sections may be read in isolation.
In the first part of each of Sections \ref{prf1}
and \ref{sectran} we give an outline of the main ideas of the
proofs before proceeding with the technical details.

\section{Preliminaries}
\label{prelim}

In this section we collect some technical results that we
need.
The first is a martingale-type
criterion for proving $\Pr [ T = \infty] >0$
for hitting times $T$.
The result is based on a well-known
idea (see e.g.~\cite[Theorem 2.2.2]{fmm} in the
countable Markov chain case).

\begin{lemma}
\label{crit2}
Let $(X_t)_{t \in \Z^+}$ be a stochastic
process on $\RR$ adapted to a filtration
$(\F_t)_{t \in \Z^+}$.
Suppose
$g: \RR \to [0,\infty)$ and
$\A \subset \RR$ (possibly infinite)
are such that
\begin{align}
\label{superm}
  \Exp [ g(X_{t+1}) - g(X_t) \mid \F_t ] \leq 0 ~{\rm on}~ \{ X_t \in \RR \setminus \A \},
  \end{align}
  for all $t \in \Z^+$.
 Write $g_0 := \inf_{x \in \A} g(x)$.
 Then for $x_0 \in \RR \setminus \A$, on $\{ X_0 = x_0 \}$
\[
\Pr [ \min \{t \in \Z^+ : X_t \in \A  \} = \infty \mid \F_0 ] \geq 1 - \frac{g(x_0)}{g_0} .\]
\end{lemma}
\proof Let $T = \min \{t \in \Z^+ : X_t \in \A  \}$,
an $(\F_t)_{t\in \Z^+}$-stopping time;
we need to show $\Pr [ T= \infty] > 0$. By (\ref{superm}),
we have that $g (X_{t \wedge T})$ is a supermartingale
adapted to $(\F_t)_{t \in \Z^+}$. Moreover, $g$ is nonnegative
so $g (X_{t \wedge T})$ converges a.s.~to some limit, say $L$.
Then if $X_0 = x_0$   we have
\[    g (x_0) \geq \Exp [ L ] \geq \Exp [ L \1_{\{T < \infty \}} ]
\geq g_0 \Pr [ T < \infty ] ,\]
which implies that $\Pr [ T < \infty ] \leq g(x_0) /g_0$, as required. \qed\\

The following maximal inequality is Lemma 3.1 in \cite{mvw}.

\begin{lemma}
\label{lm3.2}
Let  $(Y_t)_{t \in \Z^+}$
be a stochastic process on $[0,\infty)$
adapted to a filtration $(\F_t)_{t \in \Z^+}$. Suppose that $Y_0
= y_0$ and
for some $b \in (0,\infty)$ and all $t \in \Z^+$
\[  \Exp [ Y_{t+1}- Y_t \mid \F_t ] \leq b ~{\rm a.s.}. \]
Then for any $x >0$ and any $t \in \N$
\[  \Pr \left[ \max_{0 \leq s \leq t} Y_s \geq x \right]
\leq (b t +y_0) x^{-1}.    \]
\end{lemma}

Next we prove Proposition \ref{propA1}. The proof
is elementary and we do not give all the formal details.\\

\noindent
{\bf Proof of Proposition \ref{propA1}.} Clearly (A1) implies (A1$^\prime$).
So suppose that (A1$^\prime$) holds. The following four-step argument
shows that (A1) follows.

(i) First we show that without loss of generality we may take $n_{-i}
= n_i$ for each $i$. This is straightforward, since with positive probability
the walk sequentially takes $n_{-i}$ `jumps' of size $k_i$ in the $\be_i$ direction
and also with positive probability the walk sequentially takes $n_i$ `jumps'
of size $k_{-i}$ in the $-\be_i$ direction. In either case,
the walk has moved a positive distance in time $n_i n_{-i}$.

(ii) Next we show that we may take $k_i = k_{-i}$ for each $i$.
Fix $i$. In view
of part (i), we may take $n_i = n_{-i} = n$, say. Let
\[ s_i := (k_i + k_{-i}) \max \{ k_i , k_{-i} \} .\]
Without loss of generality, we may suppose $k_{-i} \geq k_i$.
Then each of the following two events has positive probability:
(a) the walk can perform $(k_i + k_{-i}) k_{-i}$
successive `jumps' of size $k_i$ in the $\be_i$ direction;
(b)  the walk can perform $(k_i + k_{-i}) k_{i}$
successive `jumps' of size $k_{-i}$ in the $-\be_i$ direction
followed by $(k_i + k_{-i}) (k_{-i} - k_{i})$
successive `jumps' of size $k_{i}$ in the $\be_i$ direction.
In either case (a) or (b), the walk ends up at distance
$(k_i + k_{-i}) k_i k_{-i}$ from its starting point
after time $n (k_i + k_{-i}) k_{-i}$.

(iii) Next we show that we can take $n_i = n_{-i} = n$ for all $i$.
Given parts (i) and (ii), we may take $k_i = k_{-i}$ and
$n_i = n_{-i}$ for each $i$. Set $n := \prod_i n_i$. For any $i$,
the walk has positive probability of performing
in succession $n/n_i$ `jumps' of size $k_i$ in either
of the $\pm \be_i$ directions. Such an event takes a total
time $n$ and leads to a positive displacement, equal in opposite directions.

(iv) Finally we show that we may take $k_i = k_{-i} = k$ for all $i$.
Given parts (i)--(iii) we can take $k_i = k_{-i}$ and $n_i = n$ for all
$i$. Set $r :=   \prod_i k_i$. Then for any $i$,
with positive probability the walk can perform
$2r/k_i$ `jumps' of size $k_i$ in the direction
$\pm \be_i$, taking time $2n r /k_i$. Then in time
$2n (r/k_i) ( k_i - 1 )$ (an even multiple of $n$)
the walk can go back and forth to achieve an additional net displacement
of $0$. The walk is then at distance $2r$ from its starting point
after a total time $2nr$. Thus (A1) holds with $n_0 = 2nr$ and $k = 2r$.
This completes the proof.
\qed\\

Finally for this section, we give the proof of Theorem \ref{thm4}.\\

 \noindent
 {\bf Proof of Theorem \ref{thm4}.} Suppose $\beta =1$.
 Suppose we are in case (a), so that $\mu(\bx) = c \| \bx \|^{-1} \be_1$.
 Then Theorem \ref{thm1} applies and $\Pr [\tau_\alpha < \infty ] =1$ for any
 $\alpha \leq \pi/2$. Now suppose we are in case (b), so that
 $\mu(\bx) = c x_1^{-1} \be_1$.
 In this case we have for any $\alpha \in (0,\pi/2)$
 \[   0 < \cos \alpha \leq  \inf_{\bx \in \W_d ( \be_1 ; \alpha )} \frac{x_1}{\| \bx \|} \leq \sup_{\bx \in \W_d ( \be_1 ; \alpha )}  \frac{x_1}{\| \bx \|} \leq 1, \]
 so that (\ref{drift1}) holds throughout $\W_d ( \be_1 ; \alpha )$. It follows from
 Theorem \ref{thm1} that $\Pr [\tau_\alpha < \infty ] =1$ for
 $\alpha < \pi/2$. Finally consider $\tau_{\pi/2}$. Let $X_t = \xi_t \cdot \be_1$; then
 $\tau_{\pi/2} = \min \{ t \in \Z^+ : X_t \leq 0 \}$. From our conditions on $\Xi$ in this case we have
 \[ \Exp [ X_{t+1} - X_t \mid \xi_t = (x_1, x_2) ] = c x_1^{-1} , ~~~ \sup_{\bx \in \W_d (\be_1 ; \pi/2 ) } \Exp [ ( X_{t+1} - X_t)^{2+\eps} \mid \xi_t = \bx ] < \infty .\]
Thus we can apply results of Lamperti \cite[Theorem 3.2]{lamp1}
to $X_t$ to conclude that $\Pr [\tau_{\pi/2} = \infty ] > 0$
for $c > c_0$ where
\[ c_0 = \frac{1}{2}  \sup_{\bx \in \W_d (\be_1 ; \pi/2 ) } \Exp [ ( X_{t+1} - X_t)^{2} \mid \xi_t = \bx ] \in (0,\infty) .\]
This completes the proof.
 \qed

\section{Finite exit times: proof of Theorem \ref{thm1}}
\label{prf1}

\subsection{Outline of the proof}

We show in this section that Theorem \ref{thm1} holds:
under the conditions of the theorem,
with probability $1$ the random walk $\Xi$ will leave
any cone $\W_d(\bu ; \alpha)$, $\alpha \in (0,\pi)$ after a finite
time. There are several steps to the proof but the
overall scheme is based on some intuitive ideas, which we now sketch.

The basic element to the proof of Theorem \ref{thm1}
is Lemma \ref{exitlem}, which says that, roughly speaking, starting
 in any small cone there is positive probability, {\em uniform
 in the current position of the walk}, that the walk hits a neighbouring
 small cone. To prove this result we need to study hitting-time
 properties of the walk. Specifically, we need to show that there
 is a good probability that the walk hits a reasonably-sized
 set at distance of the order of $\| \bx\|$ starting from $\bx$.
 The conditions (A1), (A2), and (\ref{drift1}) are of course crucial here.

In view of (A1) it is natural to work with the `$n_0$-skeleton'
process $(\xi_{t n_0})_{t \in \Z^+}$, which we denote $\Xi^*$.
In Section \ref{sec:decomp} we define a decomposition of the walk $\Xi^*$ based on the
regularity condition (A1). The basic idea is that since, by (A1),
every jump of the walk $\Xi^*$ has positive probability
of being one of $\pm k \be_i$, we can extract a symmetric random
walk from $\Xi^*$, leaving a residual process that retains
some of the regularity
of the original walk, despite no longer being Markovian.

Next, in Section \ref{hitting}, we prove our basic hitting-time
estimates.
The idea now is to treat the two parts of the decomposition
separately. The symmetric process is more straightforward to
study, and is the part of the walk that will ensure
that there is good probability of the walk hitting a particular
set some distance away without returning too close to the origin.
The technical estimate here is Lemma \ref{lm4.1}.

The next step is to show that the residual process, which has
inherited appropriate drift conditions from $\Xi$,
will with good probability not travel too far in the same time,
so that the walk as a whole has good probability of hitting
the desired set. There are complications introduced here
as the residual process depends on the realization of the
symmetric process; thus we condition on that in our estimates.
Under suitable behaviour of both processes, the walk
stays far enough from the origin that the drift remains
controlled. The technical estimate here is Lemma \ref{lm4.2}.

Based on the estimates for the two parts of our decomposition,
we show (in  Lemma \ref{lm4.3})
that $\Xi$ hits a suitable set
 with positive probability. In Section \ref{exit} we translate this
 result into our exit-from-cones result, Lemma \ref{exitlem}, which we  use
 to complete the proof of Theorem \ref{thm1}.
 Our three conditions (A1), (A2), and (\ref{drift1}) all appear
 very naturally in this scheme. Having outlined the idea, we now
 proceed with the technical work.

\subsection{Decomposition}
\label{sec:decomp}

In view of condition (A1), it is convenient to
 consider the random walk at time spacing $n_0$, i.e.~the
embedded (`skeleton')
process $(\xi_{t n_0})_{t \in \Z^+}$.
For notational convenience, set
\[ \xi^*_t := \xi_{t n_0} , ~~~(t \in \Z^+).\]
Then $\Xi^* = (\xi^*_t)_{t \in \Z^+}$
is a Markovian random walk on $\SS \subseteq \Z^d$
with transition probabilities
\[ \Pr [ \xi^*_{t+1} = \by \mid \xi^*_t = \bx ] =
\Pr [ \xi_{n_0} = \by \mid \xi_0 = \bx ] , \]
and $\xi^*_0 = \xi_0$.
 The walk $\Xi^*$
 inherits regularity from $\Xi$, as the next result shows.

 \begin{lemma}
 \label{xistar}
 (i) If (A1) holds then
  \begin{align}
  \label{minvar}
 \min_{\bx \in \SS} \min_{i \in \{1,\ldots, d\}}
   \Exp [ | ( \xi^*_{t+1} - \xi^*_t ) \cdot \be_i |^2 \mid \xi^*_t = \bx ]
  \geq 2 \kappa k^2 >0 . \end{align}
  (ii) If (A2) holds then
 \begin{align}
  \label{0508b}
  \max_{\bx \in \SS}
  \Exp [ \| \xi^*_{t+1} - \xi^*_t \|^2 \mid \xi^*_t = \bx ] \leq n_0^2 B_0 < \infty.
  \end{align}
  (iii) If (A2) and (\ref{drift1}) hold, then
  \begin{equation}
\label{stardrift2}
 \| \Exp [ \xi^*_{t+1} - \xi^*_t \mid \xi^*_t = \bx ] \| = O(\| \bx \|^{-1} ) .\end{equation}
  \end{lemma}
  \proof
  By time-homogeneity, it suffices to simplify notation by taking $t=0$ throughout.
Part (i) is immediate from (A1). For part (ii), we have by
 the triangle inequality that $\Exp [ \| \xi^*_{1} - \xi^*_0 \|^2 \mid \xi^*_0 = \bx ]$ is equal to
\begin{align*}
&   \Exp \left[ \left\| \sum_{j=1}^{n_0} ( \xi_{j} - \xi_{j-1} ) \right\|^2 \mid \xi_0 = \bx \right]
\leq \Exp \left[ \left( \sum_{j=1}^{n_0}  \| \xi_{j} - \xi_{j-1} \| \right)^2 \mid \xi_0 = \bx \right] \\
& ~~ = \sum_{j=1}^{n_0} \sum_{k=1}^{n_0}
\Exp \left[ \| \xi_j - \xi_{j-1} \| \| \xi_k - \xi_{k-1} \| \mid \xi_0 =\bx \right]
\leq \left( \sum_{j=1}^{n_0} ( \Exp [ \| \xi_j - \xi_{j-1} \|^2 \mid \xi_0 = \bx ] )^{1/2} \right)^2 ,
\end{align*}
by the Cauchy--Schwarz inequality.
  Here for $j \geq 1$, by the Markov property,
  \[ \Exp [ \| \xi_{j} - \xi_{j-1} \|^2 \mid \xi_0 = \bx ]
   = \sum_{\by \in \SS} \Exp [ \| \xi_{j} - \xi_{j-1} \|^2 \mid \xi_{j-1} = \by ]
   \Pr [ \xi_{j-1} = \by \mid \xi_0 = \bx ]
   \leq B_0 ,\]
   by (A2). Thus we obtain (\ref{0508b}).
   Finally we prove part (iii). First
   we show that   the event
   \[ E(\bx) := \left\{ \max_{0 \leq s \leq n_0}
   \| \xi_s - \bx \| > \frac{1}{2} \| \bx \| \right\} \]
   has small probability given $\xi_{0} = \bx$.
   For $\bx \in \R^d$ and $t \in \Z^+$, set $W^\bx_s  := \| \xi_s - \bx \|^2$.
   We have
   \begin{align*} \Exp [ W^\bx_{s+1} - W^\bx_s \mid \xi_s = \by ] & = \Exp [ \| \xi_{s+1} - \xi_s \|^2 + 2 (\by - \bx) \cdot ( \xi_{s+1} - \xi_s ) \mid \xi_s = \by ] \\
 &  \leq B_0 + 2 \| \by - \bx \| \| \mu ( \by) \| , \end{align*}
   using (A2). Now using (\ref{drift1}) we see that there exists $C \in (0,\infty)$ such that
   \begin{equation}
   \label{wdrift}
     \Exp [ W^\bx_{s+1} - W^\bx_s \mid \xi_s = \by ]  \leq C ( 1 + \| \bx \|^{-1} \| \by \| ) .\end{equation}
 Now define $T(\bx ) := \min \{ s \in \Z^+ : W_s^\bx > \| \bx \|^2/4 \}$.
 Note that for $s < T(\bx)$, $\| \bx \|/2 \leq \| \xi_s \| \leq 3 \| \bx \|/2$.
  Then
 given $\xi_0 = \bx$, on $\{ s < T(\bx) \}$ we have from (\ref{wdrift})
 that $\Exp [ W^\bx_{s+1} - W^\bx_{s} \mid \xi_s ] \leq 5C/2 < \infty$.
Hence we can apply Lemma \ref{lm3.2} to $W^\bx_{s \wedge T(\bx)}$ to obtain
   \begin{equation}
  \nonumber
    \Pr \left[ \max_{0 \leq s \leq n_0}
   W_{s \wedge T(\bx)}^\bx > \frac{1}{4} \| \bx \|^2 \mid \xi_{0} = \bx \right]
   \leq C \| \bx \|^{-2} ,\end{equation}
   for some $C \in (0,\infty)$. However
   $T(\bx) \leq n_0$ implies that
   \[ \max_{0 \leq s \leq n_0}
   W_{s \wedge T(\bx)}^\bx \geq W_{T(\bx)}^\bx > \frac{1}{4} \| \bx \|^2 ,\]
   by definition of $T(\bx)$. Hence
   \begin{equation}
    \label{Eprob}
     \Pr [ E( \bx ) \mid \xi_0 = \bx ] = \Pr [ T(\bx) \leq n_0 \mid \xi_0 = \bx ] \leq C \| \bx \|^{-2}.\end{equation}
    Now by   partitioning on $E(\bx)$ and applying the triangle inequality,
   \begin{align}
   \label{eqw1}
    \| \Exp [ \xi^*_{1} - \xi^*_0 \mid \xi^*_0 = \bx ]  \|
  & \leq \sum_{s=1}^{n_0} \| \Exp [ \xi_{s} - \xi_{s-1} \mid E^c ( \bx ) , \xi_{0} =\bx ] \|
     \nonumber\\
 &~~   + \sum_{s=1}^{n_0} \| \Exp [ ( \xi_{s} - \xi_{s-1} ) \1_{ E  ( \bx )} \mid \xi_{0} =\bx ] \|
,\end{align}
  where $E^c ( \bx )$ is the complementary event to $E( \bx )$.
Now using the elementary inequality that for ${\bf X}$ a random $d$-vector $\| \Exp [ {\bf X} ] \|
\leq d \Exp \| {\bf X} \|$, we have that
\[ \| \Exp [ ( \xi_{s} - \xi_{s-1} ) \1_{E ( \bx )} \mid \xi_{0} =\bx ] \|
\leq d \Exp [ \| \xi_{s} - \xi_{s-1} \| \1_{E ( \bx )} \mid \xi_{0} =\bx ] , \]
so that
by Cauchy--Schwarz, (A2), and (\ref{Eprob}),
\begin{align*}  & \| \Exp [ ( \xi_{s} - \xi_{s-1} ) \1_{E  ( \bx) } \mid \xi_{0} =\bx ] \| \\
& ~~ \leq d ( \Exp [ \| \xi_{s} - \xi_{s-1} \|^2   \mid \xi_{0} =\bx ] )^{1/2}
( \Pr [ E  (\bx ) \mid \xi_{0} = \bx  ] )^{1/2}
\leq C \| \bx \|^{-1} ,\end{align*}
for some $C \in (0,\infty)$.
On the other hand, on $E^c (\bx) \cap \{ \xi_{0} = \bx\}$
we have from (\ref{drift1}) that
\[ \max_{0 \leq s \leq n_0} \| \Exp [ \xi_{s} - \xi_{s-1} \mid E^c ( \bx ) , \xi_{0} =\bx ] \| \leq
\max_{ \by \in \SS : (1/2) \| \bx \|  \leq \| \by \| \leq (3/2) \| \bx \|} \| \mu (\by) \|^{-1} ,\]
which is again $O ( \| \bx \|^{-1})$.
Combining the two  estimates for the terms
on the right-hand side of (\ref{eqw1}) we obtain (\ref{stardrift2}).
   \qed\\

The next result establishes the decomposition. Specifically,
 we decompose the jump of $\Xi^*$ at time $t$ into a symmetric component
   ($V_{t+1}$),
 and a residual component ($\zeta_{t+1}$),  such that
 at any time $t$  only one of the two components is present in a particular
 realization.

\begin{lemma}
\label{decolem}
Suppose (A1) holds. There exist sequences of random
variables $(V_t)_{t\in \N}$ and $(\zeta_t)_{t \in \N}$
such that:
\begin{itemize}
\item[(i)]
  the $(V_t)_{t \in \N}$ are i.i.d.~with
 $V_1 \in \{\0, \pm k \be_1, \ldots, \pm k \be_d \}$ and
 \[ \Pr[V_1 = \0] = 1-2d \kappa , ~~~
 \Pr [ V_1 = -k \be_i ] = \Pr [ V_1 = +k \be_i ] = \kappa ; \]
 \item[(ii)]
  $\zeta_{t+1} \in \Z^d$ with
 $\Pr [ \zeta_{t+1} =  \0 \mid V_t \neq \0 ] =1$  and
  \begin{equation}
   \label{0512a}
   \zeta_{t+1} = \left( \xi^*_{t+1} - \xi^*_t - V_{t+1}  \right) \1_{\{V_{t+1} = \0\}}
 = ( \xi^*_{t+1} - \xi^*_t ) \1_{\{ V_{t+1} = \0\}} ; \end{equation}
 \item[(iii)]
 we can decompose
 the jumps of $\Xi^*$ via
\begin{equation}
\label{decomp}
 \xi^*_{t+1} - \xi^*_t =
       V_{t+1}   +  \zeta_{t+1}  =
   V_{t+1}  \1_{\{ V_{t+1} \neq \0\}} +  \zeta_{t+1}  \1_{\{ V_{t+1} = \0\}} ~~~(t \in \Z^+).
\end{equation}
\end{itemize}
 \end{lemma}
 \proof
 The statement of the lemma follows directly from (A1), but for clarity
let us give an explicit construction of the variables $V_t$ and $\zeta_t$.
By the time-homogeneity and Markov assumptions on $\Xi$ (hence $\Xi^*$)
for each $\bx \in \Z^d$ there exists a sequence of i.i.d.~random
vectors $\Delta^\bx_1, \Delta^\bx_2, \ldots$, independent
for each $\bx$, such that we can realize $\xi^*_{t+1} - \xi^*_t$ as
\begin{equation}
\label{deco1}
 \xi^*_{t+1} - \xi^*_t = \sum_{\bx \in \SS} \Delta^\bx_{t+1} \1_{\{ \xi^*_t = \bx \}} .\end{equation}
Condition (A1) implies that
\[ \min_{\bx \in \SS}   \Pr [ \Delta^\bx_{t+1} = k \be_i ] \geq \kappa ,\]
and similarly for $- k \be_i$. It follows that we can write
$\Delta^\bx_{t+1} = V_{t+1} + \zeta^\bx_{t+1}$ where
$V_{t+1}$ is as described in part (i) of the lemma. Then (\ref{deco1})
becomes
\[ \xi^*_{t+1} - \xi^*_t = V_{t+1} + \sum_{\bx \in \SS} \zeta^\bx_{t+1}  \1_{ \{ \xi^*_t = \bx \}} ;\]
this final sum we denote by $\zeta_{t+1}$, and parts (ii)
and (iii)
of the lemma follow. \qed\\

For $t \in \N$, (\ref{decomp})
 yields  a decomposition for $\xi^*_t$ as
\begin{equation}
\label{xidec}
 \xi^*_t - \xi_0 = \sum_{s=1}^t ( V_s   + \zeta_s ) .\end{equation}
 Note that the residual increments
 $\zeta_1, \zeta_2, \ldots$
have a rather complicated structure (and are certainly not independent); however, they do inherit
regularity properties from $\Xi$, as we summarize in the next lemma.

\begin{lemma}
\label{zetalem}
Suppose (A1) and (A2) hold, and $\zeta_{t+1}$ is as in Lemma \ref{decolem}.
Then
 \begin{align}
  \label{zetamoms}
  \max_{\bx \in \SS}
  \Exp [ \| \zeta_{t+1} \|^2 \mid \xi^*_t = \bx ] & \leq n_0^2 B_0 <\infty \\
  \label{zetamean}
  \| \Exp [ \zeta_{t+1} \mid \xi^*_t = \bx ] \|
 & =  \| \Exp [ \xi^*_{t+1} - \xi^*_t \mid \xi^*_t = \bx ] \| .
  \end{align}
\end{lemma}
\proof
From (\ref{0512a})  we have
$\| \zeta_{t+1} \| \leq \| \xi^*_{t+1} - \xi^*_t \|$ a.s., while
 conditioning on $\xi^*_t = \bx$, taking expectations on
 both sides of the first equality in (\ref{decomp})
and noting that $\Exp [ V_{t+1} \mid \xi_t^* = \bx] =\0$,
we have \[
 \Exp [ \zeta_{t+1} \mid \xi^*_t = \bx ]
= \Exp [ \xi^*_{t+1} - \xi^*_t \mid \xi^*_t = \bx ] . \]
These two facts combined with Lemma \ref{xistar} yield
the stated results. \qed

\subsection{Hitting-time estimates}
\label{hitting}

 Having established our decomposition, we will eventually use it to show
 that under the conditions of Theorem \ref{thm1}, $\Xi$ will exit
 any cone in any particular direction with good probability: see Section \ref{exit}
 below. In order to establish this result, the main ingredient will be
 the somewhat more specific Lemma \ref{lm4.3} below, which says
 that, under appropriate conditions, $\Xi$ hits some suitable ball
 with positive probability. In order to prove Lemma \ref{lm4.3}, we
 need to work separately on the two parts of the decomposition.
 We deal with the residual process in Lemma \ref{lm4.2}. First
 we study the symmetric process, building up to Lemma \ref{lm4.1}.

Set $Y_0 := \0$ and for $t \in \N$
\begin{equation}
\label{ytdef}
 Y_t := Y_0 + \sum_{s=1}^{t}   V_s . \end{equation}
 The process $( Y_t
)_{t \in \Z^+}$ is a symmetric, homogeneous
 random walk on $\Z^d$ with
$\Pr [ Y_t = Y_{t-1} ] =\Pr [ V_t =\0] = 1 - 2d\kappa < 1$ and
jumps of size $k$. For $r>0$, $\by \in \R^d$
write $B_r(\by)$ for the Euclidean $d$-ball $B_r(\by) := \{ \bx \in \R^d : \| \bx -\by \| < r \}$;
set $B_r := B_r (\0)$.
Let $\Lambda \subset \R^d$ be a convex set and $r>0$.
Define stopping times
\begin{equation}
\label{0508e}
  \sigma ( \Lambda) := \min  \{ t\in \Z^+ :  Y_t \in \Lambda  \}, ~~~
  \rho (r) :=  \min \{ t \in \Z^+ : \| Y_t \| \geq r \} . \end{equation}
  Our first  result is to show that with positive probability
  (uniformly in $N$)
  in time $\eps N^2$, for small enough $\eps$,
   the symmetric walk
    $Y_t$ hits a subset of $B_N$ of volume $\lambda N^d$ before $B_{h_0 N}$, where
    $\lambda >0$ is fixed and $h_0$ depends only on the parameters in condition (A1).

\begin{lemma}
\label{lm4.1}
Let $d \geq 2$.
Let $N \geq 1$. Let $\Lambda_N \subseteq B_N$ be a convex set
with $d$-dimensional volume at least $\lambda N^d$ for some $\lambda >0$.
Then with $\sigma$ and $\rho$ as defined at (\ref{0508e}),
there exist constants (not depending on $N$)
$N_1 \geq 1$, $h_0 \in [ 2^{-1/2} ,\infty)$, and $\eps_1 \in (0,1)$ such that
for all $N \geq N_0$, any $h \geq h_0$, and any $\eps \in (0,\eps_1)$
\[ \Pr [ \sigma ( \Lambda_N ) \leq \lfloor \eps N^2 \rfloor \wedge \rho ( h N) ] \geq \delta \]
for some $\delta >0$  depending only on $\eps$, $k$, $\kappa$, $\lambda$ and $d$, but not on $N$.
\end{lemma}
\proof Let $\eps>0$.
Note that $\Pr [ \sigma (\Lambda_N) \leq \lfloor \eps N^2 \rfloor ]
\geq \Pr [ Y_{\lfloor \eps N^2 \rfloor} \in \Lambda_N ]$.
By the standard multivariate central limit theorem for sums
of i.i.d.~random vectors, and the fact that
$\Exp [ (V_1 \cdot \be_1)^2 ] = 2 k^2 \kappa$ by Lemma
\ref{decolem}(i), we have that for measurable $A \subset \R^d$
\[ \Pr [ (2 k^2 \kappa t)^{-1/2} Y_t \in A ] \to (2 \pi)^{-d/2} \int_A \exp \{ - \| \bx \|^2/2 \} \ud \bx ,\]
as $t \to \infty$. Taking $A = (2 k^2 \kappa t)^{-1/2} \Lambda_N$ and $t = \lfloor \eps N^2 \rfloor$,
we have that the volume of $A$ is at least
$(2 k^2 \kappa \eps)^{-1/2} \lambda$, so that for some $N_1 \geq 1$ and
all $N \geq N_1$
\begin{align}
\label{clt}
 \Pr [ Y_{\lfloor \eps N^2 \rfloor} \in \Lambda_N ]
\geq \frac{1}{2} (2 \pi)^{-d/2} (2 k^2 \kappa \eps)^{-1/2} \lambda \exp \left\{ - \frac{1}{2} (2 k^2 \kappa \eps N^2)^{-1}
 \sup_{\bx \in \Lambda_N}
\| \bx \|^2 \right\} \nonumber\\
\geq \frac{1}{2} (2 \pi)^{-d/2}  (2 k^2 \kappa \eps)^{-1/2} \lambda \exp \left\{ - \frac{1}{4 k^2 \kappa \eps} \right\} ,\end{align}
since $\Lambda_N \subseteq B_N$.
 On the other hand, we claim that for any $r>0$ and $t \in \N$,
\begin{equation}
\label{maxbound}
 \Pr [ \rho ( r ) \leq t ] = \Pr \left[ \max_{0 \leq s \leq t} \| Y_s \| \geq r \right] \leq
4 d \exp \left\{ - \frac{r^2}{2d k^2 t } \right \} .\end{equation}
To obtain the  inequality in (\ref{maxbound}),
note that
\[ \Pr \left[ \max_{0 \leq s \leq t} \| Y_s \| \geq r \right] \leq d \Pr \left[ \max_{0 \leq s \leq t} | Y_s \cdot
\be_1 |   \geq d^{-1/2} r \right] ,\]
and then combine inequalities of L\'evy (see e.g.~\cite[p.~139]{gut})
and Hoeffding (see e.g.~\cite[p.~120]{gut}) on sums of i.i.d.~mean-zero bounded
random variables to obtain
\[ \Pr \left[ \max_{0 \leq s \leq t} | Y_s \cdot
\be_1 |   \geq d^{-1/2} r \right] \leq 2 \Pr \left[   | Y_t \cdot
\be_1 |   \geq d^{-1/2} r \right] \leq 4 \exp \left\{ - \frac{r^2}{2d k^2 t} \right\} .\]
Hence combining (\ref{clt}) and the $r=hN$, $t= \lfloor \eps N^2 \rfloor$ case of (\ref{maxbound})
\begin{align*}
& \Pr [ \{ \rho ( h N ) \leq \lfloor \eps N^2 \rfloor \} \cup \{ \sigma (\Lambda_N ) > \lfloor \eps N^2 \rfloor \} ] \\
& ~~ \leq
1 - \frac{1}{2} (2 \pi)^{-d/2}  (2 k^2 \kappa \eps)^{-1/2} \lambda \exp \left\{ - \frac{1}{4 k^2 \kappa \eps} \right\}
 +  4 d \exp \left \{ - \frac{h^2}{2 d k^2 \eps} \right\} \leq 1 -\delta ,\end{align*}
 for some $\delta > 0$, not depending on $N$, if we choose $h \geq h_0 := (d/\kappa)^{1/2} \geq 2^{-1/2}$
  and $\eps > 0$ small enough.
 The statement of the lemma follows.
  \qed\\

Let $Z_0 : = \0$ and for
$t \in \N$ let
\begin{equation}
\label{ztdef}
 Z_t := Z_0 + \sum_{s=1}^{t}  \zeta_s .\end{equation}
Thus $(Z_t)_{t \in \Z^+}$
is the residual part of the process $(\xi^*_t)_{t \in \Z^+}$
after the symmetric   process $(Y_t)_{t \in \Z^+}$
has been extracted. Indeed, with $Y_t, Z_t$ as defined
at (\ref{ytdef}), (\ref{ztdef}) we have
 from (\ref{xidec}) that for $t \in \N$
\begin{equation}
\label{decomp2}
 \xi^*_t -\xi_0  = Y_t   + Z_t .\end{equation}

We next show that with good probability the residual process
$(Z_t)_{t \in \Z^+}$ does not exit from a suitable ball around $\0$
by time $\lfloor \eps N^2\rfloor$. By construction
the process $(Z_t)_{t \in \Z^+}$ depends upon $  ( V_t) _{t \in \N}$
because   the distribution of $\zeta_{t+1}$ depends upon the
value of $\xi^*_{t}$. For $t \in \N$,
 let $\Omega_V(t) :=
  \{\0, \pm k \be_i ,\ldots, \pm k \be_d \}^{t}$ and let $\omega_V
\in \Omega_V(t)$ denote a generic
realization of the sequence  $( V_1, \ldots, V_t)$.
For $r >0$ define
\[ \Omega_{V,r} ( t ) = \{ \omega_V \in \Omega_{V}(t) : t < \rho (r ) (\omega_V)  \} , \]
i.e., the set of those paths $\omega_V$ for which $\| Y_s \| < r$
for all $s \leq t$.
Our next result, Lemma \ref{lm4.2}, gives control over the deviations
of $Z_t$.
The choice of $3/4$ as the lower bound
in Lemma \ref{lm4.2}  is fairly arbitrary: any lower bound in
$(0,1)$ can be obtained for $\eps$ small enough, but $3/4$ is good enough for us.

\begin{lemma}
\label{lm4.2}
Let $h \in (0,\infty)$. Suppose   (A1)
and (A2) hold and that for some $K_0 \in (0,\infty)$
\begin{equation}
\label{mudisc}
 \max_{\bx \in \SS \cap B_{(1+h) N} ( \xi_0 )} \| \Exp [ \xi^*_{t+1} - \xi^*_t \mid \xi^*_t =\bx \| \leq K_0 N^{-1} ,
\end{equation}
for all $N \geq 1$.
Let $c \in (0,1/2]$.
There exists $\eps_2 > 0$ not depending
on $N$ (but depending on $c$  and $K_0$) such that
for all $\eps \in (0,\eps_2)$, all $N \geq 1$,
and all $\omega_{V} \in \Omega_{V,h N} ( \lfloor \eps N^2 \rfloor)$
\[ \Pr \left[ \max_{0 \leq t \leq \lfloor \eps N^2 \rfloor}
\| Z_t  \| \leq
c  N \mid (V_1,\ldots, V_{\lfloor \eps N^2 \rfloor} ) = \omega_V, Z_0 = \0    \right] \geq \frac{3}{4} .
\]
\end{lemma}
\proof
Let $c \in (0,1/2]$.
 For the duration of this proof, define the stopping time
\begin{align*}
\tau_0 := \min \{ t \in \Z^+ : \| Z_t  \| > c N \}.
\end{align*}
For the duration of this proof,
let $\G_t = \sigma (\xi^*_0, \xi^*_1, \ldots, \xi^*_t , V_1, \ldots, V_t)$.
Then $\zeta_1, \ldots, \zeta_t$ and $Z_0, Z_1, \ldots, Z_t$
are $\G_t$-measurable, and $\tau_0$ is a $(\G_t)_{t \in \Z^+}$ stopping time.
Consider
the stopped square-deviation
process defined for $t \in \Z^+$ by
$W_t := \| Z_{t \wedge \tau_0 }   \|^2$; $W_t$ is then $\G_t$-adapted.
Suppose that $t \leq \lfloor \eps N^2\rfloor$.
On the event $\{ \tau_0 > t \}$ we have that
\begin{align*}
W_{t+1} - W_t = \| Z_{t +1}   \|^2
-  \| Z_{t }  \|^2 =
\| Z_{t+1}-Z_t \|^2 + 2  ( Z_{t+1} - Z_t ) \cdot
 Z_t
= \| \zeta_{t+1} \|^2 + 2 \zeta_{t+1} \cdot Z_t  ,\end{align*}
while on $\{ \tau_0 \leq t\}$, $W_{t+1} - W_t = 0$.
So
conditioning on $\G_t$
and taking expectations,
 we obtain
\begin{align}
\label{0512b}
  \Exp [ W_{t+1} - W_t \mid   \G_t ]
=   \left( \Exp [ \| \zeta_{t+1} \|^2 \mid   \G_t ]
+ 2 \Exp [  \zeta_{t+1} \cdot Z_t \mid  \G_t ] \right) \1_{\{ t < \tau_0 \}} .\end{align}
The first term on the right-hand side
of (\ref{0512b}) is at most
\[  \Exp [ \| \zeta_{t+1} \|^2 \mid   \G_t ]
= \Exp [ \| \zeta_{t+1} \|^2 \mid \xi^*_t , V_{t+1} ] =
O(1), \]
 by
 (\ref{zetamoms}).
 For the second  term on the right-hand side of
 (\ref{0512b}), since $Z_t$ is a measurable function of $\G_t$,
  \begin{align}
 \label{0806a}
| \Exp [ \zeta_{t+1} \cdot Z_t  \mid  \G_t ] | \1_{\{ t < \tau_0 \}}
 & \leq  \| Z_t \| \left\| \Exp [ \zeta_{t+1} \mid \G_t ] \right\|    \1_{\{ t < \tau_0 \}}   \nonumber\\
 & \leq  c N  \left\| \Exp [ \zeta_{t+1} \mid \G_t ] \right\|    \1_{\{ t < \tau_0 \}} , \end{align}
 by the definition of $\tau_0$.
 Now we have
 \begin{align*} \| \Exp [ \zeta_{t+1} \mid \G_t ]
 \1_{\{ t < \tau_0 \}}
  \| & \leq \esssup_{A \in \G_t : \, t < \tau_0 (A)}
 \| \Exp [ \zeta_{t+1} \mid \xi^*_t = \xi^*_t(A), V_{t+1} = V_{t+1} (A)  ] \| ,\\
 & \leq \esssup_{A \in \G_t : \, t < \tau_0 (A)}
 \| \Exp [ \zeta_{t+1} \mid \xi^*_t = \xi^*_t(A), V_{t+1} = \0  ] \| , \end{align*}
 since, by (\ref{0512a}), $\zeta_{t+1} = \zeta_{t+1} \1_{\{V_{t+1} = \0\} }$.
 By the same fact,  \begin{align*}
 \Exp [ \zeta_{t+1} \mid V_{t+1} = \0, \xi^*_t ]
 = (  \Pr[ V_{t+1} = \0])^{-1} \Exp [ \zeta_{t+1} \mid   \xi^*_t ]
  = (1-2d\kappa)^{-1}  \| \Exp [ \zeta_{t+1}
  \mid \xi^*_t ] \|
  , \end{align*}
   by Lemma \ref{decolem}(i). Combining the last two displayed equations, we have
   that there exists $C = C( d, \kappa) < \infty$ such that
   \[ \| \Exp [ \zeta_{t+1} \mid  \G_t ] \| \1_{\{ t < \tau_0  \}} \leq
   C \esssup_{A \in \G_t : \, t < \tau_0 (A)  } \| \Exp [ \zeta_{t+1}
  \mid \xi^*_t = \xi^*_t (A) ] \| .\]
  Now suppose $\omega_V \in \Omega_{V,hN} ( \lfloor \eps N^2 \rfloor )$. Then
  from (\ref{decomp2})
  \[ \max_{0 \leq s \leq t} \| \xi^*_s - \xi_0 \| \leq \max_{0 \leq s \leq t} \| Z_s \|
  + \max_{0 \leq s \leq t} \| Y_s \| \leq c N + h N \leq (1+h) N \]
  on $\{ (V_1,\ldots, V_{\lfloor \eps N^2 \rfloor} ) = \omega_V ,
  t < \tau_0 \}$. In particular, using (\ref{zetamean}), this implies that
  \[ \esssup_{A \in \G_t : \, t < \tau_0(A), \, \omega_V(A) = \omega_V } \| \Exp [ \zeta_{t+1}
  \mid \xi^*_t = \xi^*_t (A) ] \| \leq \sup_{\bx \in B_{(1+h)N} (\xi_0)  }  \| \Exp [ \xi^*_{t+1} -\xi^*_t
  \mid \xi^*_t = \bx ] \| .\]
   Hence, assuming (\ref{mudisc}), we obtain, for any $\omega_V \in \Omega_{V,hN} ( \lfloor \eps N^2 \rfloor )$,
   \begin{align}
 \label{0806b}
  \| \Exp [  \zeta_{t+1}  \mid \omega_V, \G_t ] \| \1_{\{ t < \tau_0 \}} \leq
     C K_0 N^{-1} .\end{align}
  Thus combining (\ref{0806a})
 and (\ref{0806b}) we have
 \begin{equation}
 \label{0512c}
|   \Exp [  \zeta_{t+1} \cdot Z_t \mid \omega_V, \G_t ] |
\1_{\{ t < \tau_0  \}}
 \leq d^{1/2} c N C K_0 N^{-1} = C d^{1/2} c     K_0   .\end{equation}
 Hence from (\ref{0512b}) with (\ref{0512c}) and (\ref{zetamoms}) we have, a.s.,
\[ \sup_{\omega_V \in \Omega_{V , hN} ( \lfloor \eps N^2 \rfloor )}
\Exp [ W_{t+1} - W_t \mid   \omega_V, \G_t  ] \leq    C d^{1/2} c K_0 + n_0^2 B_0 =: B_1 ,\]
where $B_1 \in (0,\infty)$ does not depend on $\eps$ or $N$.
 Then
applying Lemma \ref{lm3.2}   we have
\begin{align*}
 \Pr\left[ \max_{0 \leq t \leq \lfloor \eps N^2 \rfloor} W_t \geq c^2 N^2
  \mid \omega_V \right]
   \leq \frac{B_1 \eps N^2}{c^2 N^2} = c^{-2} \eps B_1.
   \end{align*}
So taking $\eps_2 = c^2/(4B_1)$ and $\eps \in (0,\eps_2)$, we have
\begin{equation}
\label{0510a}
\Pr\left[ \max_{0 \leq t \leq \lfloor \eps N^2 \rfloor } \| Z_{t \wedge \tau_0 }  \| \leq c N
  \mid \omega_V \right] =
 \Pr\left[ \max_{0 \leq t \leq \lfloor \eps N^2 \rfloor} W_t \leq c^2 N^2
  \mid \omega_V \right]    \geq \frac{3}{4} .\end{equation}
But since, by definition of $\tau_0$,
$\| Z_{\tau_0}   \| > cN$,
we have that the left-hand
event in (\ref{0510a})
implies that $\tau_0   > \lfloor \eps N^2 \rfloor$, and so
we obtain the required result.
\qed\\

 Lemmas \ref{lm4.1} and \ref{lm4.2}
 give us control over the two parts
 of the decomposition of $\Xi^*$.
 Our final ingredient
  before we can prove
  the main result of this section (Lemma \ref{lm4.3})
  is the next lemma, which gives control over the deviations
 of $\Xi$ from the embedded process $\Xi^*$.

 \begin{lemma}
 \label{xidev}
  Suppose that (A2) holds.
 There exist $\eps_3 > 0$ and $N_2 \geq 1$  such that
for all $\eps \in (0,\eps_3)$, all $N \geq N_2$,
and all $\omega_{V} \in \Omega_{V} ( \lfloor \eps N^2 \rfloor)$
 \[ \Pr \left[ \max_{0 \leq s \leq n_0 \lfloor \eps N^2 \rfloor } \| \xi_s - \xi_{ \lfloor s / n_0 \rfloor} \|
 \leq \frac{N}{8} \mid (V_1,\ldots, V_{\lfloor \eps N^2 \rfloor} ) = \omega_V \right] \geq \frac{3}{4} .\]
 \end{lemma}
 \proof
 We
have
\[ \max_{0 \leq s \leq n_0 \lfloor \eps N^2 \rfloor  } \| \xi_s - \xi_{ \lfloor s / n_0 \rfloor} \|
\leq \max_{0 \leq r \leq \lfloor \eps N^2 \rfloor } \max_{1 \leq s \leq n_0 - 1}
\| \xi_{n_0r +s} - \xi_{n_0 r} \| ,\]
where, by the triangle inequality,
\[ \max_{1 \leq s \leq n_0 - 1}
\| \xi_{n_0r +s} - \xi_{n_0 r} \|
= \max_{1 \leq s \leq n_0 - 1}
\left\| \sum_{j=0}^{s-1} \xi_{n_0r +j+1} - \xi_{n_0 r +j } \right\|
\leq \sum_{j=0}^{n_0 -2} \| \xi_{n_0r +j+1} - \xi_{n_0 r +j }  \| .\]
Thus to complete the proof of the lemma, we need to show that
\begin{equation}
\label{1111}
 \Pr \left[ \max_{0 \leq r \leq \lfloor \eps N^2 \rfloor } \sum_{j=0}^{n_0 - 2} \| \xi_{n_0 r + j+1}
-\xi_{n_0+r +j} \|   >  \frac{N}{8} \mid \omega_V \right]  < \frac{1}{4}, \end{equation}
for suitable $\eps$, $N$ and all $\omega_V$. For each $r$ we have, by Cauchy--Schwarz,
\[ \Exp \left[ \left( \sum_{j=0}^{n_0 - 2} \| \xi_{n_0 r + j+1}
-\xi_{n_0+r +j} \| \right)^2 \mid \omega_V \right] \leq
\left( \sum_{j=0}^{n_0 - 2}   ( \Exp [ \| \xi_{n_0 r +j+1} - \xi_{n_0 r +j} \|^2 \mid \omega_V ] )^{1/2}
\right)^2  , \]
and the expectation here satisfies
\[ \Exp [ \| \xi_{n_0 r +j+1} - \xi_{n_0 r +j} \|^2 \mid \omega_V ]
= \sum_{\bx \in \SS} \Exp [ \| \xi_{n_0 r +j+1} - \xi_{n_0 r +j} \|^2 \mid \xi_{n_0 r + j} =\bx ]
\Pr [ \xi_{n_0 r + j} =\bx \mid \omega_V ] \leq B_0 ,\]
by (A2). Hence, by Boole's inequality followed by Markov's inequality,
 the probability on the left-hand side of (\ref{1111}) is bounded above by
\[ \sum_{0 \leq r \leq \lfloor \eps N^2 \rfloor }
\Pr \left[ \sum_{j=0}^{n_0 - 2} \| \xi_{n_0 r + j+1}
-\xi_{n_0+r +j} \|   >  \frac{N}{8} \mid \omega_V \right] \leq 64 N^{-2} (1 + \lfloor \eps N^2 \rfloor )
n_0^2 B_0 .\]
This is less than $1/4$ for $\eps < 2^{-9} n_0^{-2} B_0^{-1}$ and $N \geq 2^5 n_0 B_0$; thus we verify (\ref{1111})
and the lemma follows.
 \qed\\

Combining the preceding three lemmas,
we can prove the key result of this section, Lemma \ref{lm4.3}, which says
 that with positive probability $\Xi$ hits a
  sizable $d$-ball
 in $B_N(\xi_0)$ before it leaves the ball $B_{2h_0 N}(\xi_0)$; this is the next result.

\begin{lemma}
\label{lm4.3}
Let $d \geq 2$, $N \geq 1$, and $\xi_0 \in \SS$. Suppose
that (A1) and (A2) hold and that for some $K_0 \in (0,\infty)$ (\ref{mudisc})
holds with $h = 2h_0$, where $h_0$ is the constant in Lemma \ref{lm4.1}.
   Let $c \in (0,1)$ and
   $\Lambda_N = B_{cN/4} (\xi_0+\by) \subseteq B_N(\xi_0)$, for some $\by \in B_{3N/4}$.
Then
there exists $\delta >0$ and $N_0 \geq 1$ (neither depending on $N$)
such that for all $N \geq N_0$
\[ \Pr [ \Xi \textrm{ hits } \Lambda_N \textrm{ before exit from } B_{2 h_0 N} (\xi_0)  ] \geq \delta .\]
\end{lemma}
\proof
Let $\eps_0 = \eps_1 \wedge \eps_2 \wedge \eps_3$
and $N_0 = \max \{ N_1, N_2 \}$. Take $\eps \in (0,\eps_0)$ and $N \geq N_0$.
Let $c \in (0,1)$. Let $h_0 \geq 2^{-1/2}$
be as in Lemma \ref{lm4.1}. Fix $\by \in B_{3N/4}$
so that $\Lambda_N = B_{cN/4} (\xi_0+\by) \subseteq B_N(\xi_0)$. Also
let $\Lambda'_N = B_{cN/8} (\xi_0+\by) \subseteq B_N(\xi_0)$.
Define the events
\begin{align*}
 G := \left\{ \max_{0 \leq t \leq \lfloor \eps N^2 \rfloor} \| Z_t  \| \leq c N/8 \right\}, ~~~
H := \left\{ \sigma ( \Lambda'_N ) \leq \lfloor \eps N^2 \rfloor \wedge \rho ( h_0 N) \right\}, \\
{\rm and} ~~~ I:= \left\{ \max_{0 \leq s \leq n_0 \lfloor \eps N^2 \rfloor } \| \xi_s - \xi_{ \lfloor s / n_0 \rfloor} \|
 \leq \frac{N}{8} \right\} .\end{align*}
 Write $\sigma = \sigma (\Lambda_N')$.
 Then on $H$, $\sigma \leq \lfloor \eps N^2 \rfloor$ and
 $\| Y_{\sigma} - \by \| \1_{H} \leq c N/8$ so that
 on $G \cap H$
\[ \| \xi^*_{\sigma} - \xi_0 - \by \| \leq \| Y_{\sigma} - \by \| + \| Z_{\sigma} \| \leq (c N/8) + (c N/8) = c N/ 4.\]
Thus $\xi^*_{\sigma} = \xi_{n_0 \sigma} \in \Lambda_N$ on $G \cap H$.
Next we need to control $\| \xi_s - \xi_0 \|$ for $s$ up to $n_0 \lfloor \eps N^2 \rfloor$.
For any $t$ we have
\begin{equation}
\label{max2}
 \max_{0 \leq s \leq n_0 t } \| \xi_s - \xi_0 \|
\leq \max_{0 \leq s \leq n_0 t} \| \xi_s - \xi_{ \lfloor s / n_0 \rfloor} \|
+ \max_{0 \leq s \leq n_0 t} \|  \xi_{ \lfloor s / n_0 \rfloor} - \xi_0 \|  .\end{equation}
For $t= \lfloor \eps N^2 \rfloor$, the first term on the right-hand side of (\ref{max2})
is bounded by $N/8$ on $I$.
For the second term on the right-hand side of (\ref{max2}),
it follows from (\ref{decomp2}) and the
triangle inequality
that
\[
\max_{0 \leq s \leq n_0 \lfloor \eps N^2 \rfloor } \|  \xi_{ \lfloor s / n_0 \rfloor} - \xi_0 \|  \1_{ G \cap H}
\leq
\max_{0 \leq t \leq \lfloor \eps N^2 \rfloor }
\| \xi^*_t - \xi_0 \| \1_{ G \cap H} \leq h_0 N + (c N/8) .\]
Thus, from (\ref{max2}),
\[  \max_{0 \leq s \leq n_0 t } \| \xi_s - \xi_0 \| \1_{G \cap H \cap I} \leq (N/8) + h_0 N + (c N/8) \leq 2 h_0 N ,\]
since $h_0 \geq 2^{-1/2}$.
Hence (with $\xi_0$ as given)
\[ E: = \left\{  \Xi \textrm{ hits } \Lambda_N \textrm{ before exit from } B_{2 h_0 N} (\xi_0)
\right\} \supseteq G \cap H \cap I.\]
$H$ is determined by the realization
$\omega_V \in \Omega_V (\lfloor \eps N^2 \rfloor)$, and  so (with $\xi_0$ as given)
\begin{align*} \Pr [ E   ] \geq
\Pr[ G \cap H \cap I  ]
= \sum_{\omega_V \in \Omega_V (\lfloor \eps N^2 \rfloor ) : H ~{\rm occurs}}
  \Pr[ G \cap I
\mid \omega_V  ] \Pr[  \omega_V    ] \\
= \sum_{\omega_V \in \Omega_{V,h_0 N} (\lfloor \eps N^2 \rfloor ) : H ~{\rm occurs}}
  \Pr[ G \cap I
\mid \omega_V  ] \Pr[  \omega_V    ] ,\end{align*}
by definition of $H$ and $\Omega_{V,h_0 N} (\lfloor \eps N^2 \rfloor )$.
But from Lemma \ref{lm4.2} with $c =1/8$ and Lemma \ref{xidev} we have that $\Pr [ G \cap I \mid \omega_V ] \geq 1/2$
for all $\omega_V \in \Omega_{V,h_0N} (\lfloor \eps N^2 \rfloor )$, since $\eps < \eps_2 \wedge \eps_3$
and $N \geq N_2$. Hence we obtain
 \begin{align*}
 \Pr[ E  ] \geq
\frac{1}{2} \sum_{\omega_V \in \Omega_{V,h_0N} (\lfloor \eps N^2 \rfloor ) : H ~{\rm occurs} }
  \Pr[    \omega_V    ]
  = \frac{1}{2} \Pr [ H   ] \geq \delta/2 > 0 ,
  \end{align*}
 applying Lemma \ref{lm4.1}, since $\eps < \eps_1$ and $N \geq N_1$.
     \qed\\

     \rem
     At first glance, one might hope to prove Lemma \ref{lm4.3}
 by choosing $\eps$ small enough in Lemmas \ref{lm4.2} and \ref{xidev}
 so that we can replace the lower bounds of $3/4$ there by something
 very close to $1$, and then combine this with Lemma \ref{lm4.1} to show
 that $G \cap H \cap I$ (as in the proof above) occurs with positive probability
 using a simple union bound. This does not work, however, since as $\eps$
 gets small, the $\delta$ in Lemma \ref{lm4.1} gets smaller too. That is why
 we needed to use the more sophisticated argument, conditioning on the
 path of $Y_t$.

 \subsection{Exit from cones}
 \label{exit}

The next result is essentially a restatement of Lemma \ref{lm4.3}
 in the context in which we will apply it
 to complete
 the proof of Theorem \ref{thm1}.

\begin{lemma}
\label{exitlem}
Let $d \geq 2$. Suppose
that (A1) and (A2) hold. Suppose (\ref{drift1}) holds.
Then for any $c \in (0,1)$,
there exist  $A_1 \in (0,\infty)$ and $\delta >0$
such that
\[ \min_{\bx \in \SS : \| \bx \| \geq A_1} \min_{\by \in \SS :  \| \by - \bx \| \leq (a_0/2) \| \bx \|}
\Pr [ \Xi \textrm{ hits } B_{(ca_0/6) \| \bx \|} (\by ) \mid \xi_0 = \bx ] \geq \delta ,\]
where $a_0 \in (0,1)$ is a constant that does not depend on $c$.
\end{lemma}
  \proof
  Suppose $\xi_0 = \bx \in \SS$. Take $h = 2h_0$, where $h_0$ is the constant in Lemma \ref{lm4.1}.
  Set $N = \frac{1}{2(1+h)} \| \bx \|$ and take $\| \bx \|$ large enough so that $N \geq 1$.
  Now assuming (\ref{drift1}), we have from (\ref{stardrift2}) that for some $C \in (0,\infty)$
  \[ \| \Exp [ \xi^*_{t+1} - \xi^*_t \mid \xi^*_t = \by ] \| \leq C \| \by \|^{-1} ,\]
  for all $\by \in \SS$, so that, since $(1+h)N = \| \bx \|/2$,
  \[ \sup_{\by \in \SS \cap B_{(1+h)N} (\bx) } \| \Exp [\xi^*_{t+1} - \xi^*_t \mid \xi^*_t = \by ] \|
  \leq 2 C \| \bx \|^{-1} ,\]
  uniformly in $\bx$. In other words, (\ref{mudisc}) holds for some
  $K_0$ and all $N \geq 1$. Take $a_0 = \frac{3}{4(1+h)} < 1$. Then
  $\| \by - \bx \| \leq a_0 \| \bx \|/2$
  implies that
   $\| \by - \bx \| \leq 3N/4$. Hence setting
   $\Lambda_N = B_{(c a_0 /6) \| \bx \|} (\by)$,
   where $\by \in B_{3N/4} (\bx)$ and $(ca_0/6) \| \bx \| = N/4$,
    Lemma \ref{lm4.3} is applicable; therefore the result follows
   for all $N \geq N_0$, that is, for
   $\| \bx \| \geq 2 (1+2h_0) N_0 = A_1$, say.
  \qed\\

Now we can complete the proof of Theorem \ref{thm1}.\\

\noindent
{\bf Proof of Theorem \ref{thm1}.} We show that for arbitrary $\bu \in \S_d$ and arbitrary $\eps>0$, $\Xi$ eventually hits
$\W_d(\bu; \eps)$ in finite time with probability 1. Without loss of generality, fix $\eps>0$ (small) and consider the
cone $\W_d (\be_1 ; \eps)$: we want to show that eventually $\Xi$ enters this cone.
Given $\eps$, with $a_0$ the constant in Lemma \ref{exitlem},
take $c =(4/a_0) \tan \eps \in (0,1)$, for $\eps$ small enough.
Then let $A_1$ be the constant given by Lemma \ref{exitlem} with this choice of $c$.
For any $d$ and any $\eps$, we can find a finite set $\{ \bu_1 = \be_1, \bu_2, \ldots, \bu_k \} \subset \S_d$
and $\eps' \in (0,\eps)$ such that
\[ B_{A_1}  \cup \left(
\bigcup_{j=1}^k \W_d ( \bu_j ; \eps ) \right) = \R^d , \]
but where
\[ \W_d ( \bu_i ; \eps') \cap \W_d ( \bu_j ; \eps' ) \cap \{ \bx : \| \bx \| > A_1 \} = \emptyset \]
 for all $i \neq j$.
 Denote $C_0 := B_{A_1} $ and for $i \in \{1, \ldots, k \}$, $C_i := \W_d ( \bu_i ; \eps) \setminus C_0$,
 $C'_i := \W_d ( \bu_i ; \eps') \setminus C_0$.
 If $C_i \cap C_j \neq \emptyset$ for $i \neq j$, we say that
 $i$ and $j$ are neighbours. For neighbours $i$ and $j$,
 we have (for small enough $\eps$) that in the notation of Lemma \ref{exitlem},
 with $c$ small enough, for any $\bx \in C_i$ we can always find $\by$ with $\| \by - \bx \| \leq (a_0/2) \| \bx\|$
 such that $B_{(c a_0 /6) \| \bx \|} (\by) \subset C'_j$. Hence an application of Lemma
 \ref{exitlem} yields that for neighbours $i$ and $j$
 \begin{equation}
 \label{irred} \Pr [ \Xi \textrm{ hits $C_j$ } \mid \xi_t \in C_i ] \geq \delta >0 ,\end{equation}
 where $\delta$ does not depend on $i$, $j$, or $\xi_t$.

 Define a $\{0,1,\ldots,k\}$-valued stochastic process
 $(J_t)_{t \in \Z^+}$ by $J_t := \min \{ j : \xi_t \in C_j \}$.
 Condition (A1) ensures that if $J_t = 0$ then with positive
 probability $J_r >0$ for some $r >t$. Moreover, (\ref{irred})
  implies that uniformly in the location of $\xi_t$,
 there is positive probability that after time $t$
 $\Xi$ hits a neighbouring
 cone of $C_{J_t}$. The state-space of $J_t$ is finite,
 and by the above argument state $0$ is not absorbing
 while all the non-zero states communicate. It follows by standard `irreducibility'
 arguments
 that $J_t$ hits any non-zero state in finite time with probability $1$,
 and in particular $J_T =1$ for some $T < \infty$. This completes the proof. \qed

\section{Limiting direction: proof of Theorem \ref{thm2}}
\label{sectran}

\subsection{Overview and notation}

 The aim of this section
is to prove Theorem \ref{thm2}, and demonstrate
the existence of a limiting direction.
We will deduce Theorem \ref{thm2} from the following result on exit from cones
for the walk $\Xi$, which says that under the conditions of Theorem \ref{thm2},
provided $\Xi$ starts `far enough inside'
a cone, there is probability close to $1$ that it remains in the cone for all time.

\begin{theorem}
\label{thm6}
Let $d \in \{2,3,\ldots\}$ and $\bu \in \S_d$.  Suppose that (A2) holds and that
for some $\beta \in (0,1)$, $c >0$, $\delta >0$, and $A_0 >0$,
(\ref{trancona}) and (\ref{trancon2a}) hold.
Let $\alpha \in (0,\pi)$ and $\eps>0$. Then there exists $\alpha' \in (0,\alpha)$ (not depending on $\eps$)
and
$A_1 < \infty$ such that for any $\bx \in \SS \cap \W_d ( \bu ; \alpha')$ with
$\| \bx \| > A_1$
\[ \Pr [ \tau_\alpha = \infty \mid \xi_0 = \bx ] \geq 1 - \eps .\]
\end{theorem}

The scheme for the
proof of Theorem \ref{thm6} is as follows.
First, we prove a two-dimensional version of Theorem \ref{thm6}, that says
 for any two-dimensional cone (`wedge'),
under suitable conditions,
$\Pr [ \tau_\alpha = \infty ] \geq 1-\eps$. To prove Theorem \ref{thm6} on exit from cones in general $d \geq 2$,
we use an argument
based on projections
down onto two-dimensional subspaces.  In order to apply the projection argument,
we need to extend the two-dimensional walks that we consider
from  Markov processes to processes that are adapted
to some larger filtration. Thus now we
establish the relevant formalism, and then state our two-dimensional
result, Theorem \ref{thm5}.

 For $\bx = (x_1,x_2) \in \R^2$
we use the notation $\bx_\perp = (-x_2, x_1)$.
Let $(\F_t)_{t \in \Z^+}$ be a filtration. Suppose that $Z = (Z_t)_{t \in \Z^+}$
is an $(\F_t)_{t \in \Z^+}$-adapted process on $\R^2$.
For what follows, we will typically take $\F_t$ to be $\sigma (\xi_1,\ldots,\xi_t)$
for the random walk $\Xi$ on $\Z^d$ and take $Z_t$ to be an appropriate
 projection onto $\R^2$
of $\xi_t$.
For our results on $Z$,   we assume
 the following regularity condition
analogous to (A2).
\begin{itemize}
 \item[(A3)] There exists $B_0 \in (0,\infty)$ such that
 \[ \max_{t \in \Z^+}
 \esssup \Exp [ \| Z_{t+1} - Z_t \|^{2} \mid \F_t ] \leq B_0, \]
 where the essential supremum is over all $A \in \F_t$ with $\Pr (A) >0$.
 \end{itemize}

We are now ready to state the two-dimensional result that will allow
us to deduce Theorem \ref{thm6} and hence Theorem \ref{thm2}.

 \begin{theorem}
 \label{thm5}
 Let $d=2$.
 Suppose that (A3) holds.
 Let $\alpha \in (0,\pi)$ and $\bu \in \S_d$.
 Suppose
that for some $\beta \in (0,1)$, $c>0$, $\delta>0$,
   $A_0>0$, and some $(\F_t)_{t \in \Z^+}$-stopping time $\sigma$,
 \begin{align}
 \label{trancon}
  \min_{\bx \in \SS \cap \W_2( \bu; \alpha ) : \| \bx \| > A_0 } \min_{t \in \Z^+}
 \essinf_{\{ Z_t = \bx \} \cap \{ t < \sigma \}} \left(
 \| \bx \|^{\beta} \Exp [ Z_{t+1} - Z_t \mid \F_t ] \cdot \hat \bx \right)  \geq c ,
 ~\textrm{and} \\
 \label{trancon2}
 \max_{\bx \in \SS \cap \W_2 (\bu ; \alpha) : \| \bx \| > A_0 } \max_{t \in \Z^+}
 \esssup_{\{ Z_t = \bx \} \cap \{ t < \sigma \}}
 \left(
  \| \bx \|^{\beta+\delta} | \Exp [ Z_{t+1} - Z_t \mid \F_t ] \cdot \hat \bx_\perp | \right) < \infty .
 \end{align}
 Fix $\eps>0$.
 Then
 there exist $\alpha' \in (0,\alpha)$ and $A_1 < \infty$ such that for
 any $\bx \in \SS \cap \W_2 ( \bu ; \alpha' )$
 with $\| \bx \| > A_1$
 \[ \Pr [ \min \{ t \in \Z^+ : Z_{t \wedge \sigma} \notin \W_2 (\bu ; \alpha ) \} = \infty \mid \F_0 ] \geq 1 - \eps  \]
 on $\{ Z_0 = \bx \}$.
 \end{theorem}

 \rem We could state Theorem \ref{thm6} (and indeed Theorem \ref{thm2})
  at a similar level of generality
 as Theorem \ref{thm5}, i.e., replacing $\Xi$ with a more general
 adapted process $Z_t$ on $\Z^d$. However, this extra generality is unnecessary for the main
 line of this section, which is the proof of Theorem \ref{thm2}.

\subsection{Proof of Theorem \ref{thm5}}

In this section we prove Theorem \ref{thm5}.
For the moment we restrict our attention
to the problem of exit from the
quadrant $Q := \W_2(\be_1 ; \pi/4) = \{ (x_1, x_2) \in \R^2 : x_1 >0, |x_2| < x_1 \}$, where the computations
are more transparent.
It will be convenient to use polar coordinates for $\bx = (x_1,x_2)$,
so that $x_1 = r \cos \varphi$,
$x_2 = r \sin \varphi$ where $r = \| \bx \|$ and $\varphi$ is the angle between $\bx$ and
$\be_1$, measured anticlockwise.
For  $\nu>0$ and $\bx \in \R^2$
set
\begin{align}
\label{hdef}
h_{\nu} (\bx) = h_{\nu} (r,\varphi ) :=
r^{-2\nu} (\cos (2\varphi) )^{-1}
= \frac{(x_1^2 +x_2^2)^{1-\nu}}{x_1^2-x_2^2}.\end{align}
Then $h_{\nu}$ is positive in the interior of the quadrant $Q$ and
 blows up on the boundary $\partial Q$.

\begin{figure}
\centering
\includegraphics[width=10cm]{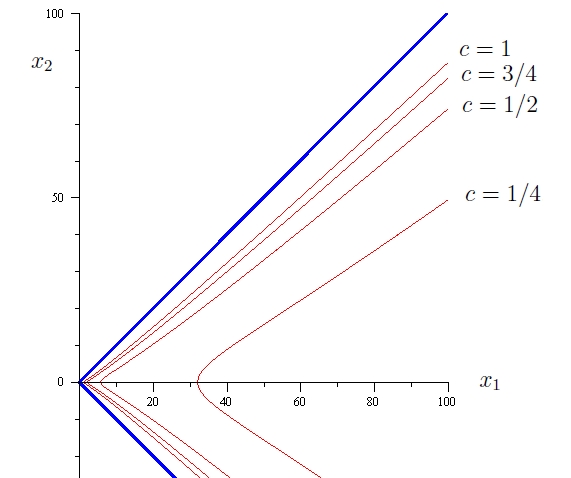}
\caption{Plot of segments of the contours
 $\gamma_{0.2} (c)$ for $4c \in \{1,2,3,4\}$. The $c=1/4$
 contour cuts the $x_1$-axis at $(1/4)^{-1/0.4} = 32$. }
\label{contours}
\end{figure}

For $s>0$ define the unbounded
open subset of $\R^2$
\[ \Gamma_{\nu} (s) := \{\bx \in \R^2  : 0 < h_{\nu} (\bx) < s\} .\]
Then for $t \geq s > 0$,
  $\Gamma_{\nu} (s) \subseteq \Gamma_{\nu} (t) \subseteq
  Q$,
 and for any $s >0$,
$x_1 \to\infty$
 as $\| \bx \| \to \infty$
along any path in $\Gamma_\nu(s)$. Note that the contours \[
\gamma_\nu (c) :=
\{\bx \in Q : h_\nu (\bx) = c\} = \partial \Gamma_\nu (c) ,\]
$c>0$, eventually leave any wedge $\W_2 ( \be_1; \beta )$, $\beta \in (0, \pi/4)$, and
so approach the boundary of $Q$ in this angular sense. However, they do so relatively slowly.
In particular, an elementary calculation shows that for fixed $\nu$
and fixed $c_1 > c_2 >0$, for $\bx \in \gamma_\nu (c_1)$
\begin{equation}
\label{hasymp}
 \inf_{\by \in \gamma_\nu (c_2) }
\| \bx - \by \| \sim ( c_2^{-1} - c_1^{-1} ) \| \bx \|^{1-2\nu} ,\end{equation}
as $\| \bx \| \to \infty$, so that the contours diverge the farther out into the wedge
they go. Also observe that $\gamma_\nu (c)$
cuts the $x_1$-axis at $(c^{-1/(2\nu)},0)$.
See Figure \ref{contours} for an example.

 Given $\bx \in \Gamma_{\nu} (s)$ we have
 from (\ref{hdef}) that
 \begin{align}
 \label{hh1}
 \| \bx\|^{-2\nu} \leq h_\nu ( \bx) \leq s   .
\end{align}
We work with a truncated version of
 $h_\nu$, namely $\tilde h_{\nu} : \R^2 \to [0,1]$,
 defined for $\bx \in \R^2$ by
 \[  \tilde h_{\nu} (\bx) := \begin{cases}
 \min \{ h_\nu (\bx) , 1 \} & ~~\textrm{for}~~ \bx \in Q; \\
 1 & ~~\textrm{for}~~ \bx \in \R^2 \setminus Q. \end{cases}
 \]
 Observe that for $s>0$
 \begin{equation}
 \label{hinf}
  \inf_{\bx \notin \Gamma_\nu (s) } \tilde h_\nu (\bx) = \min \{ 1, s\} .\end{equation}

 We will derive some basic properties of the functions $h_{\nu}$
 and $\tilde h_\nu$. To this end, we will use multi-index notation for partial derivatives on $\R^2$. For
$\sigma = (\sigma_1,\sigma_2) \in \Z^+ \times \Z^+$, $D_\sigma$ will
denote $D_1^{\sigma_1} D_2^{\sigma_2}$
where $D_j^k$ for $k \in \N$ is $k$-fold
differentiation with respect to $x_j$, and $D^0_j$ is the identity
operator. We also use the notation $|\sigma|:=\sigma_1+\sigma_2$
and $\bx^\sigma :=  x_1^{\sigma_1} x_2^{\sigma_2}$.

\begin{lemma}
Let $\nu \in (0,1)$ and $s \in (0,1)$. Then for $\bx \in \Gamma_\nu (s)$
and $\by = (y_1,y_2)$
\begin{align}
\label{hh2}
  \sum_{j=1}^2 y_j D_j h_\nu (\bx)   =   - 2 \nu \| \bx \|^{-1} h_\nu (\bx)
  \left( ( \by \cdot \hat \bx )
  - 2 \nu^{-1} x_1 x_2 \| \bx \|^{2\nu-2} h_\nu (\bx) ( \by \cdot \hat \bx_\perp )   \right) .\end{align}
  Also there exists $C \in (0,\infty)$
  such that for any $\bx \in \Gamma_\nu (s)$ and $\by = (y_1,y_2)$
  \begin{equation}
  \label{hh2a}
   \left| \sum_{j=1}^2 y_j D_j h_\nu (\bx)  \right| \leq C \| \by \| \| \bx \|^{2\nu -1} h_\nu (\bx) .\end{equation}
Moreover
for any $\bx \in \Gamma_\nu (s)$,
as $\| \bx\| \to \infty$
\begin{equation}
\label{hh3}
\sup_{ \sigma : |\sigma|=2}
| D_\sigma   h_{\nu} (\bx) | = O (  \| \bx \|^{4\nu -2}  h_\nu (\bx)   ) .
\end{equation}
\end{lemma}
\proof
Let $\nu, s \in (0,1)$.
Directly from (\ref{hdef})
we obtain
\begin{align}
\label{hd1}
D_1   h_\nu (\bx) & = \frac{2(1-\nu) x_1 (x_1^2 +x_2^2)^{-\nu}}{x_1^2-x_2^2} -
\frac{2x_1 (x_1^2+x_2^2)^{1-\nu}}{(x_1^2-x_2^2)^2}, ~~\textrm{and}~~\nonumber\\
D_2   h_\nu (\bx) & = \frac{2(1-\nu) x_2 (x_1^2 +x_2^2)^{-\nu}}{x_1^2-x_2^2} +
\frac{2x_2 (x_1^2+x_2^2)^{1-\nu}}{(x_1^2-x_2^2)^2}.\end{align}
Since for $\bx = (r, \varphi)$ in polar coordinates, for any
$\by = (y_1,y_2)$,
\begin{align*}
 y_1   = (\by \cdot \hat \bx) \cos \varphi -
(\by \cdot \hat \bx_\perp )  \sin \varphi,
~~\textrm{and}~~
y_2   = (\by \cdot \hat \bx ) \sin \varphi +
(\by \cdot \hat \bx_\perp )  \cos \varphi,
\end{align*}
it follows from (\ref{hd1}) that
\begin{align*}
\sum_{j=1}^2 y_j D_j h_\nu (\bx)
 & =   -\frac{2 \nu (x_1^2+x_2^2)^{(1/2)-\nu}}{x_1^2-x_2^2} ( \by \cdot \hat \bx )
+ \frac{4 x_1 x_2 (x_1^2+x_2^2)^{(1/2)-\nu}}{(x_1^2-x_2^2)^2} ( \by  \cdot \hat \bx_\perp )  \\
 & =   -\frac{2 \nu (x_1^2+x_2^2)^{(1/2)-\nu}}{x_1^2-x_2^2}  \left(
 \by \cdot \hat \bx
- 2 \nu^{-1} x_1 x_2 \| \bx \|^{2\nu-2} h_\nu (\bx) \by \cdot \hat \bx_\perp    \right) ,\end{align*}
which yields (\ref{hh2}). Now from (\ref{hh2}) we have that
\begin{align*}
\left| \sum_{j=1}^2 y_j D_j h_\nu (\bx)  \right| \leq C \| \by \| \left( \| \bx \|^{-1}
h_\nu (\bx) + h_\nu (\bx)^2 \| \bx \|^{2 \nu - 1} \right) \\
\leq C \| \by \| \| \bx \|^{-1} h_\nu (\bx ) \left( 1 + h_\nu (\bx) \| \bx \|^{2 \nu} \right) ,
\end{align*}
which with (\ref{hh1}) yields (\ref{hh2a}).
Similarly, differentiating in (\ref{hd1})
and using   (\ref{hh1})
we obtain (\ref{hh3}).
\qed\\

  We next show that when
  (\ref{trancon})
  holds,
 $(\tilde h_{\nu}(Z_t))_{t \in \Z^+}$ is a
supermartingale on $\Gamma_{\nu} (s)$
 for   suitably
small $\nu, s > 0$. This is the next result.

\begin{lemma}
\label{hsup}
Suppose that (A3) holds.
Suppose that for some $\beta \in (0,1)$, $c >0$, $\delta>0$,
  $A_0 >0$, and $(\F_t)_{t \in \Z^+}$-stopping time $\sigma$,
(\ref{trancon}) and (\ref{trancon2}) hold.
Then there exist $\nu, s \in (0,1/2)$ such that
for any $t \in \Z^+$
 \[ \Exp [ \tilde h_{\nu} (Z_{t+1}) - \tilde h_{\nu} (Z_t) \mid \F_t ] \leq 0 \]
  on $\{ Z_t \in \Gamma_{\nu} (s) \} \cap \{ t < \sigma \}$.
\end{lemma}
 \proof
 We suppose throughout that $t < \sigma$.
 Let $\nu >0$ be such that $\nu < \min \{ \delta/2, (1-\beta)/8\} < 1/8$.
 Let $s \in (0,1/2)$, to be fixed later.
 Note that, by (\ref{hh1}), if $\bx \in \Gamma_{\nu} (s)$
 we have $\| \bx \| > s^{-1/(2\eps)}$ and $\tilde h_{\nu} (\bx ) = h_\nu (\bx)$.
 Also note that since $\tilde h_{\nu} (\bx) \in [0,1]$ for all $\bx$, we
 have
 \begin{equation}
 \label{qq21}
  \left| \left( \tilde h_{\nu} (\bx + \by) - \tilde h_{\nu} (\bx) \right)
 - \1_{\{ \| \by \| <   \| \bx \|^{1-3\nu} \}} \left( \tilde h_{\nu} (\bx + \by) -
 \tilde h_{\nu} (\bx) \right)
 \right|
\leq \1_{\{ \| \by \| \geq   \| \bx \|^{1-3\nu} \}} ,\end{equation}
for any $\bx, \by \in \R^2$. We have from (\ref{hasymp}) that there exists
$C_1 = C_1 (s,\nu) \in (0,\infty)$ such that for all $\bx \in \Gamma_\nu (s)$
with $\| \bx \| > C_1$, for any $\by$ with
 $\| \by \| < \| \bx \|^{1-3\nu}$,
$\bx + \by \in \Gamma_\nu (2s) \subset \Gamma_\nu (1)$. Thus
 Taylor's theorem
with    Lagrange form for the remainder
implies that for $\bx \in \Gamma_\nu(s)$ with $\| \bx \| > C_1$,
\begin{align}
\1_{\{ \| \by \| < \| \bx \|^{1-3\nu} \}}
\left( \tilde h_{\nu} (\bx + \by) -\tilde  h_{\nu} (\bx) \right) =
\1_{\{ \| \by \| < \| \bx \|^{1-3\nu} \}}
\left(  h_{\nu} (\bx + \by) -  h_{\nu} (\bx) \right) \nonumber\\
= \1_{\{ \| \by \| < \| \bx \|^{1-3\nu} \}} \sum_{j=1}^2 y_j (D_j    h_{\nu}  )(\bx)
+ \frac{1}{2} \1_{\{ \| \by \| <   \| \bx \|^{1-3\nu} \}}
 \sum_{\sigma: |\sigma|=2} \by ^\sigma (D_\sigma  h_{\nu} ) (\bx + \eta \by ),
 \label{qq22}
 \end{align}
for some $\eta = \eta (\by) \in (0,1)$. Taking $\bx = Z_t$ and $\by = Z_{t+1} - Z_t$
and combining (\ref{qq21}) and (\ref{qq22}) we have
that on $\{ Z_t \in \Gamma_\nu (s), \| Z_t \| \geq C_1 \}$,
\begin{align}
 \Exp [ \tilde h_{\nu} (Z_{t+1}) - \tilde h_{\nu} (Z_t) \mid \F_t ]
  = \Exp \left[ \1_{\{ \| \by \| < \| \bx \|^{1-3\nu} \}}
  \sum_{j=1}^2 y_j (D_j   h_{\nu}  )(\bx) \mid \F_t \right] \nonumber\\
+ \frac{1}{2}
\Exp \left[ \1_{\{ \| \by \| <   \| \bx \|^{1-3\nu} \}}
 \sum_{\sigma: |\sigma|=2} \by ^\sigma (D_\sigma    h_{\nu} ) (\bx + \eta \by ) \mid \F_t \right]
 + K \Pr [ \| \by \| \geq   \| \bx \|^{1-3\nu}  \mid \F_t ],
 \label{qq22a}
 \end{align}
 where $| K| \leq 1$.

 We now deal with each of the terms on the right-hand side
 of (\ref{qq22a}) in turn.
For the final term on the right-hand side of (\ref{qq22a}),
the conditional
form of Markov's   inequality and (A3) give, for $\bx = Z_t$ and $\by = Z_{t+1} - Z_t$,
\begin{align}
\label{jumpbound}
 \Pr [ \| \by \| \geq   \| \bx \|^{1-3\nu}  \mid \F_t ]
 \leq \| \bx \|^{6 \nu-2}   \Exp [ \| Z_{t+1} - Z_t \|^2 \mid \F_t ]
\leq B_0 \| \bx \|^{6 \nu-2} .
\end{align}
 The first term on the right-hand side of (\ref{qq22a}) may be written as
\begin{align*}
 \Exp \left[  \sum_{j=1}^2 y_j (D_j   h_{\nu}  )(\bx) \mid \F_t \right]
  - \Exp \left[ \1_{\{ \| \by \| \geq \| \bx \|^{1-3\nu} \}}
  \sum_{j=1}^2 y_j (D_j   h_{\nu}  )(\bx) \mid \F_t \right] ,\end{align*}
  where by (\ref{hh2a}) we have
  \[ \left| \Exp \left[ \1_{\{ \| \by \| \geq \| \bx \|^{1-3\nu} \}}
  \sum_{j=1}^2 y_j (D_j   h_{\nu}  )(\bx) \mid \F_t \right] \right|
  \leq C \| \bx \|^{2 \nu -1} h_\nu (\bx)
   \Exp \left[ \1_{\{ \| \by \| \geq \| \bx \|^{1-3\nu} \}} \| \by \| \mid \F_t \right] .\]
   By Cauchy--Schwarz, this last expression is bounded by
   \[ C \| \bx \|^{2 \nu -1} h_\nu (\bx) \left( \Pr [ \| \by \| \geq \| \bx \|^{1-3\nu}  \mid \F_t ]
   \right)^{1/2} ( \Exp [ \| \by \|^2 \mid \F_t ] )^{1/2} = O ( \| \bx \|^{5 \nu - 2}
   h_\nu (\bx) ) ,\]
   by (\ref{jumpbound}) and (A3). For the second term on the right-hand side
   of (\ref{qq22a}), we have from (\ref{hh3}) that
   \begin{align*}
   \left| \Exp \left[ \1_{\{ \| \by \| <   \| \bx \|^{1-3\nu} \}}
 \sum_{\sigma: |\sigma|=2} \by ^\sigma (D_\sigma    h_{\nu} ) (\bx + \eta \by ) \mid \F_t \right]
 \right| \\
 \leq C  \| \bx + \eta \by \|^{4 \nu -2} h_\nu (\bx + \eta \by ) \1_{\{ \| \by \| <   \| \bx \|^{1-3\nu} \}}
 = O(  \| \bx \|^{4\nu -2} h_\nu (\bx) ) ,\end{align*}
 for $\bx \in \Gamma_\nu (s)$ with $\| \bx \| > C_1$. Combining these
 calculations we obtain from (\ref{qq22a}) that
 \begin{equation}
 \label{tt1} \Exp [ \tilde h_{\nu} (Z_{t+1}) - \tilde h_{\nu} (Z_t) \mid \F_t ]
  = \Exp \left[  \sum_{j=1}^2 y_j (D_j   h_{\nu}  )(\bx) \mid \F_t \right]
  + O ( \| \bx \|^{6 \nu -2} ),\end{equation}
 on $\{ Z_t = \bx \}$ for $\bx \in \Gamma_\nu (s)$ with $\| \bx \| > C_1$, and where $\by = Z_{t+1} - Z_t$.

Now from  (\ref{trancon}) and (\ref{trancon2}) we have that for $\| \bx \| > A_0$,
\[ \Exp [ \by \cdot \hat \bx \mid \F_t ] \geq c \| \bx \|^{-\beta} , ~~~
 | \Exp [ \by \cdot \hat \bx_\perp \mid \F_t ] | = O ( \| \bx \|^{-\beta -\delta }) .\]
Hence taking expectations in (\ref{hh2}),   on $\{ Z_t = \bx \}$
for $\bx \in \Gamma_\nu (s)$ with $\| \bx\|$ large enough,
\begin{align}
\label{tt2}
\Exp \left[  \sum_{j=1}^2 y_j (D_j   h_{\nu}  )(\bx) \mid \F_t \right]
\leq -2 \nu \| \bx \|^{-1} h_\nu (\bx) \left[
c \| \bx \|^{-\beta} + O ( \| \bx \|^{2 \nu}   \| \bx \|^{-\beta-\delta} ) \right] \nonumber\\
\leq - 2 \nu \| \bx \|^{-1} h_\nu (\bx) \| \bx \|^{-\beta} ( c + o(1) ) ,\end{align}
since $\nu < \delta/2$.
Noting that, by (\ref{hh1}), for $\bx \in \Gamma_\nu (s)$ we can replace
the $O ( \| \bx \|^{6 \nu -2} )$ term in
(\ref{tt1}) by $O  ( \| \bx \|^{8 \nu -2} h_\nu (\bx) )$, we obtain
from (\ref{tt1}) and (\ref{tt2})
\[ \Exp [ \tilde h_{\nu} (Z_{t+1}) - \tilde h_{\nu} (Z_t) \mid \F_t ]
\leq - 2 \nu \| \bx \|^{-1} h_\nu (\bx) \| \bx \|^{-\beta} ( c + o(1) + O ( \| \bx \|^{\beta+8 \nu -1} ) ),\]
which is negative for all $\| \bx\|$ large enough, since $\nu < (\beta-1)/8$.
Also, for $\bx \in \Gamma_\nu (s)$ we have from
(\ref{hh1}) that $\| \bx \| \geq s^{-1/(2\nu)}$.
So taking
$s$ small enough,
the result follows.
\qed\\

\noindent
{\bf Proof of Theorem \ref{thm5}.}
It suffices to consider wedges with principal
axis in direction $\be_1$.
First we prove
the theorem for the quadrant case, $\alpha =\pi/4$.
In this case,
Lemma \ref{hsup}
shows that $\tilde h_{\nu}( Z_{t \wedge \sigma} )$ is a supermartingale
in $\Gamma_{\nu} (s)$
for $\nu,s$ small enough. Choose $\nu, s \in (0,1/2)$ as in Lemma \ref{hsup}, take some $K>1$
(to be fixed later) and write
$\Gamma := \Gamma_\nu (s)$, $\Gamma' := \Gamma_\nu (s/K) \subset \Gamma$ for this choice
of parameters.
Then by (\ref{hinf}) and the definition of $\Gamma'$,
\[ \inf_{\bx \in Q  \setminus \Gamma } \tilde h_{\nu} ( \bx) \geq s, ~~~{\rm and }~~~ \sup_{\bx \in \Gamma' } \tilde h_{\nu} ( \bx) \leq s/K .\]
Thus
Lemma \ref{crit2} applies with $X_t = Z_{t \wedge \sigma}$
and $g = \tilde h_\nu$.  Thus
 for any $\bx \in \Gamma'$, on $\{ Z_0 = \bx\}$,
\[
\Pr [ \min \{ t \in \Z^+ : Z_{t \wedge \sigma } \notin \Gamma  \} = \infty \mid \F_0 ] \geq 1 - \frac{s/K}{s} = 1 - \frac{1}{K} .\]
This in turn implies that for any $\bx \in \Gamma'$, on $\{ Z_0 = \bx \}$,
\begin{equation}
\label{quadtran}
  \Pr [ \min \{ t \in \Z^+ : Z_{t \wedge \sigma} \notin Q  \} = \infty  \mid \F_0 ] \geq 1 - \frac{1}{K} ,\end{equation}
which we can make as close to $1$ as we like by choosing $K$ large enough. Finally, since the contours
$\gamma_\nu(c)$ eventually leave any wedge inside $Q$, we note that
given $K, \nu, s$ and $\theta \in (0,\pi/4)$ we can find $A_1$ large enough such that
$\{ \bx \in \W_2(\be_1 ; \theta) : \| \bx \| > A_1\} \subseteq \Gamma'$.
This proves Theorem \ref{thm5} for $\alpha = \pi/4$, and hence any $\alpha \geq \pi/4$ too.

Now we extend this argument to   angles $\alpha \in (0,\pi/4)$.  For such $\alpha$,
let ${\bf L}_\alpha$
denote the linear transformation of $\R^2$ defined by
\[ {\bf L}_\alpha = \left( \begin{array}{ll} \cos \alpha & 0 \\
0 & \sin \alpha   \end{array} \right)  .\]
Then ${\bf L}_\alpha \W_2 (\be_1 ; \pi/4) = \W_2 (\be_1 ; \alpha)$.

So consider the random walk $Z_t$ in wedge
$\W_2 (\be_1 ; \alpha)$. Given that condition
(\ref{trancon}) holds in $\W_2 (\be_1 ; \alpha)$,
the same condition also holds
for the walk ${\bf L}_\alpha^{-1} ( Z_t )$ on $\W_2(\be_1 ; \pi/4)$.
  Hence the argument for (\ref{quadtran})
implies that for small enough $\nu, s$ and
for any $\bx \in {\bf L}_\alpha \Gamma_\nu (s/K)$, on $\{ Z_0 =\bx\}$,
\[
 \Pr [ \min \{ t \in \Z^+ : Z_{t \wedge \sigma} \notin \W_2 (\be_1 ; \alpha)  \} = \infty  \mid \F_0 ] \geq 1 - \frac{1}{K} ,\]
and we argue as previously.
 This completes the proof of the theorem.
 \qed

 \subsection{Proof of Theorem \ref{thm6}}
 \label{projections}

 \noindent
 {\bf Proof of Theorem \ref{thm6}.}
 The case $d=2$ of Theorem \ref{thm6} is immediate from Theorem \ref{thm5} on taking
 $Z_t = \xi_t$ and $\sigma = \infty$. So suppose
 $d \in \{3, 4,\ldots\}$.
 It suffices to work with cones with principal axis in the $\be_1$ direction and with angle
 $\alpha >0$ small (but fixed). Write $C = \W_d(\be_1 ; \alpha)$, $d>2$.
  We want to show that $\Xi$ remains in
 $C$ with probability close to $1$ if it starts far enough `inside' the cone.
 Let $\pi_1, \ldots, \pi_{d-1}$ be two-dimensional projections from $\R^d$ defined by
  $\pi_j : (x_1,\ldots,x_d) \mapsto (x_1,x_{j+1})$,
where $j \in \{1, \ldots, d-1\}$.

 For $R \subseteq \R^d$ write $\pi_j(R) \subseteq \R^2$ for its projection and
 $\Pi_j(R)$ for the inverse image $\pi_j^{-1} ( \pi_j (R)) \subseteq \R^d$, i.e.,
 $\Pi_j(R) := \{ (x_1,\ldots ,x_d) \in \R^d : (x_{1},x_{j+1}) \in \pi_j (R) \}$.
 For cones such as $C$, $\pi_j(C)$ is a wedge
 (a copy of $\W_2(\be_1 ; \alpha)$) in $\R^2$ and $\Pi_j(C)$
 is a copy of $\pi_j (C) \times \R^{d-2}$. In particular,
 $\bx \in \cap_{j=1}^{d-1} \Pi_j (C)$ implies that
 $\bx = (x_1,\ldots, x_d)$ satisfies $x_1 >0$ and
 $d-1$ linear inequalities each involving $x_1$ and one of
 $x_2, \ldots, x_d$. Thus $\cap_{j=1}^{d-1} \Pi_j (C)$ is a
 convex rectilinear cone  that contains the
 circular cone $\W_d( \be_1 ; \alpha)$. By an elementary geometrical
 argument, and convexity, the rectilinear cone $\cap_{j=1}^{d-1} \Pi_j (C)$
 is contained in a circular cone $\W_d( \be_1 ; \alpha_0)$ for some $\alpha < \alpha_0 < c(d) \alpha$,
 where $c(d)$ is a constant depending only on the dimension $d$.

 In particular, this argument shows that there exists $\alpha' \in (0,\alpha)$
 such that the
 $d$-dimensional circular
 cone $C' = \W_d(\be_1 ; \alpha') \subset C$ satisfies
 \[ C \supseteq \cap_{j=1}^{d-1} \Pi_j (C') .\]
Thus the event
\[ E:= \bigcap_{j=1}^{d-1} \{ \pi_j ( \xi_t ) \in \pi_j ( C' ) ~{\rm for~all~} t \} \]
implies that $\xi_t \in C$ for all $t$, that is, $\tau_\alpha = \infty$.
Thus it suffices to show that for any $\eps>0$ we have $\Pr [ E] \geq 1- \eps$
provided $\xi_0 \in C'' = \W_d (\be_1 ; \alpha'')$, with $\| \xi_0\|$ large enough, for some $\alpha'' \in (0,\alpha')$.
Here
\begin{equation}
\label{eq22}
 \Pr [ E ] \geq 1 - \sum_{j=1}^{d-1} \Pr [ \pi_j ( \xi_t ) ~{\rm exits~from~} \pi_j(C') ] .\end{equation}

Let $Z^{(j)}_t = \pi_j (\xi_t)$ for $j \in \{1,\ldots,d-1\}$.
Define the corresponding exit times
\[ T_j = \min \{ t \in \Z^+ : Z^{(j)}_t \notin \pi_j ( C' ) \} ,\]
so that $\cap_{j=1}^{d-1}\{ T_j = \infty \}$ implies $\{ \tau_\alpha = \infty \}$.
Given $\xi_0 \in C''$ we have that $Z^{(j)}_0 \in \pi_j (C'')$, which
is a wedge strictly contained in $\pi_j (C')$.
Thus  Theorem \ref{thm5} applies with $\sigma = \tau_\alpha$,
an $(\F_t)_{t \in \Z^+}$-stopping time. Hence
 there exist the putative $\alpha'' \in (0,\alpha')$ and $A_1$ such that
if $\| Z_0 \| > A_1$, with probability at least $1-(\eps/d)$
the process $Z_{t \wedge \tau_\alpha}$ remains inside $\pi (C')$. The same argument applies
to each of the $d-1$ probabilities in (\ref{eq22}), and so
we have
that with probability at least $1 -\eps$, $Z^{(j)}_{t \wedge \tau_\alpha} \in \pi_j (C')$
for all $t$ and all $j$. This implies that either (i)
$\tau_\alpha < \infty$ and $Z^{(j)}_{\tau_\alpha} \in \pi_j (C')$ for all $j$, or
(ii) $\tau_\alpha = \infty$. However, case (i) is impossible
since by construction $Z^{(j)}_{\tau_\alpha} \in \pi_j (C')$ for all $j$
implies that $\xi_{\tau_\alpha} \in C$, which is a contradiction by the definition
of $\tau_\alpha$. Thus we conclude that
$\Pr [ \tau_\alpha = \infty ]  \geq 1 - \eps$.
This completes the proof. \qed

 \subsection{Proof of Theorem \ref{thm2}}

 To complete the proof of Theorem \ref{thm2},
 we deduce from Theorem \ref{thm5} the existence
 of a limiting direction.\\

\noindent
{\bf Proof of Theorem \ref{thm2}.}
 Fix $\alpha>0$ (small). We show that for any $\bv \in \S_d$,
 there is positive probability that
 the walk eventually remains within angle $\alpha$ of $\bv$. Thus fix $\bv \in \S_d$.
 With this $\alpha$, let $A_1$ and $\alpha' \in (0,\alpha)$ be the
 constants in the $\eps=1/2$ case of Theorem \ref{thm6}.
 Then for some $K = K(d,\alpha') \in \N$ there exists
 a set $\{ \bu_1,   \ldots, \bu_{K-1}, \bu_K   \} \subset \S_d$
 such that
 \[ \R^d =  \bigcup_{i=1}^K \W_d ( \bu_i ; \alpha' )  ,\]
 i.e., we can write $\R^d$ as the union
 of $K$ cones (labelled $1,\ldots, K$)
 of interior half-angle $\eps/2$.
 Each of the cones $\W_d( \bu_i ; \alpha' )$
 sits inside the larger cone
 $\W_d ( \bu_i ; \alpha )$.
 Write
 $B := B_{A_1}$
 for the ball of radius $A_1$.
 We use the notation
 \[ \W_d' ( \bu ; \, \cdot \,  ) := \W_d ( \bu ; \, \cdot \, ) \setminus B .\]
 Consider the
 stochastic process
 \[ \Theta (t) := \max \{ 1 \leq j \leq K : \xi_t \in \W_d' (\bu_j ; \alpha' ) \},
 \]
 with the convention that $\max \emptyset := 0$. Thus
 $\Theta(t) = 0$ if and only if $\xi_t \in B$; otherwise
 $\Theta(t)$ takes the label of one of the truncated
 cones $\W_d' (\bu_j ; \alpha' )$
 containing $\xi_t$.

Condition (A1) implies that
\[ \min_{\bx \in \SS \cap B}
\min_j \Pr [ \Xi \textrm{ hits } \W'_d ( \bu_j ; \alpha' ) \mid \xi_0 = \bx ] \geq p > 0 .\]
It follows that a.s.~there exist infinitely many times $t_1, t_2, \ldots$
for which $\xi_{t_j} \notin B$. For each $t_j$, we have $\Theta (t_j) > 0$.
If $\Theta (t_j) >0$,
Theorem \ref{thm6}
 implies that with   probability at least $1/2$ (uniformly in
 $j$ and $\xi_{t_j}$) the walk
 remains in the larger cone $\W_d ( \bu_{\Theta (t_j)} ; \alpha )$
 for all time $t \geq t_j$.
 It follows that: (i) eventually $\Xi$ remains in some
 cone $\W_d ( \bu_{\Theta} ; \alpha )$, where $\Theta =
 \lim_{t \to \infty} \Theta (t)$; and
 (ii) $\Pr [ \Theta = j ] >0$ for all $j \in \{1 ,\ldots, K\}$.
One consequence of (i) is that $\| \xi_t \| \to \infty$ a.s.,
i.e., the walk is transient.

In other words, (i) says that, for any $\alpha>0$, eventually the walk remains within
 angle $\alpha$ of some $\bu_\Theta$, so
 $\xi_t / \| \xi_t \|$ has an
almost sure limit. Moreover, (ii) says that   with positive
probability $\xi_t / \| \xi_t \|$ remains arbitrarily
close to any of the $\bu_j$, and in particular to the given
vector $\bv$.  Thus the limit in question has
 distribution supported on all of $\S_d$. \qed

 \begin{center}
 {\bf Acknowledgements}
 \end{center}

 Some of this work was done when AW was at the
 University of Bristol, partially supported
by the Heilbronn Institute for Mathematical Research.


\begin{thebibliography}{99}

\bibitem{aim} S. Aspandiiarov, R. Iasnogorodski, and M. Menshikov,
Passage-time moments for nonnegative stochastic processes and an application to reflected random walks in a quadrant,
{\em Ann. Probab.} {\bf 24} (1996) 932--960.

\bibitem{bansmi} R. Ba\~nuelos and R.G. Smits, Brownian motion in cones,
{\em Probab. Theory Relat. Fields} {\bf 108} (1997) 299--319.

\bibitem{burkh} D.L. Burkholder, Exit times of Brownian motion,
harmonic majorization, and Hardy spaces,
{\em Adv. Math.} {\bf 26} (1977) 182--205.

 \bibitem{cohenbook} J.W. Cohen,
 Analysis of Random Walks,
 {\em Studies in Probability, Optimization and Statistics} {\bf 2},
 IOS Press, Amsterdam, 1992.

 \bibitem{dante2} D. DeBlassie and R. Smits,
 The influence of a power law drift on the exit
 time of Brownian motion from a half-line,
 {\em Stochastic Processes Appl.} {\bf 117} (2007) 629--654.

 \bibitem{fmm} G. Fayolle, V.A. Malyshev, and M.V. Menshikov,
Topics in the Constructive Theory of Countable Markov Chains,
Cambridge University Press, 1995.

\bibitem{gut} A. Gut, Probability: A Graduate Course,
Springer, 2005.

\bibitem{karlsson} A. Karlsson,
Linear rate of escape and convergence in direction,
pp.~459--471 in: Random Walks and Geometry, V.A. Kaimanovich (Ed.),
de Gruyter, 2004.

 \bibitem{klein} L.A. Klein Haneveld and A.O. Pittenger,
 Escape time for a random walk from an orthant,
 {\em Stochastic Processes Appl.} {\bf 35} (1990) 1--9.

\bibitem{lamp1} J. Lamperti,
Criteria for the recurrence and
transience of stochastic processes I,
{\em J. Math. Anal. Appl.} {\bf 1} (1960) 314--330.

\bibitem{lamp2} J. Lamperti, Criteria for
stochastic processes II: passage-time moments,
{\em J. Math. Anal. Appl.} {\bf 7} (1963) 127--145.

 \bibitem{lawler91} G.F. Lawler, Estimates for differences and Harnack
 inequality for difference operators coming from random walks with
 symmetric, spatially inhomogeneous, increments,
 {\em Proc. London Math. Soc.} {\bf 63} (1991) 552--568.

   \bibitem{lawlerbook} G.F. Lawler,
   Intersections of Random Walks,
   {\em
   Probability and Its Applications}, Birkh\"auser, Boston,
   1996.

   \bibitem{lawlerconf} G.F. Lawler,
    Conformally Invariant Processes in the Plane,
   American Mathematical Society, 2005.

   \bibitem{wedges} I.M. MacPhee, M.V. Menshikov, and A.R. Wade,
   Exit times from cones for non-homogeneuos random
   walk with asymtptoically zero drift,
   Preprint \texttt{arXiv:0806.4561}.

    \bibitem{mai} M.V. Menshikov, I.M. Asymont, and R. Iasnogorodskii,
Markov processes with asymptotically zero drifts,
{\em Problems of
Information Transmission} {\bf 31} (1995) 248--261; translated from
{\em Problemy Peredachi Informatsii} {\bf 31} (1995) 60--75 (in
Russian).

 \bibitem{mp} M.V. Menshikov and S.Yu. Popov,
 Exact power estimates for countable Markov chains,
 {\em Markov Processes Relat. Fields} {\bf 1} (1995) 57--78.

 \bibitem{mvw} M.V. Menshikov, M. Vachkovskaia, and A.R. Wade,
  Asymptotic behaviour of randomly reflecting billiards
  in unbounded tubular domains,
 {\em J. Stat. Phys.} {\bf 132} (2008) 1097--1133.

  \bibitem{mw1}
  M.V. Menshikov and A.R. Wade,
  Random walk in random environment with asymptotically zero perturbation,
{\em J. Euro. Math. Soc.} {\bf 8} (2006) 491--513.

  \bibitem{mustapha} S. Mustapha,
  Gaussian estimates for spatially inhomogeneous random walks
  on $Z^d$,
  {\em Ann. Probab.} {\bf 34} (2006) 264--283.

   \bibitem{shi} Z. Shi,
   Windings of Brownian motion and random walks in the plane,
   {\em Ann. Probab.} {\bf 26} (1998) 112--131.

  \bibitem{spitzer} F. Spitzer, Some theorems concerning $2$-dimensional
  Brownian motion,
  {\em Trans. Amer. Math. Soc.} {\bf 87} (1958) 187--197.

  \bibitem{spitzerrw} F. Spitzer, Principles of Random Walk,
  2nd edition, Springer, New York, 1976.

\bibitem{var1} N.Th. Varopoulos, Potential theory in conical domains,
{\em Math. Proc. Camb. Phil. Soc.} {\bf 125} (1999) 335--384.

\bibitem{zeitouni} O. Zeitouni,
Random walks in random environments,
{\em J. Phys. A} {\bf 39} (2006) R433--R464.

\end{thebibliography}
\end{document}